\documentclass[12pt]{amsart}
\usepackage{fullpage}
\usepackage{bm}

\usepackage{amsfonts}
\usepackage{amssymb}
\usepackage{amscd}
\usepackage{color}
\usepackage{graphicx}
\usepackage{graphics}
\usepackage{hyperref}
\input xy
\xyoption{all}

\newtheorem{theorem}{Theorem}[section]
\newtheorem{lemma}[theorem]{Lemma}
\newtheorem{corollary}[theorem]{Corollary}
\newtheorem{proposition}[theorem]{Proposition}

\theoremstyle{definition}
\newtheorem{definition}[theorem]{Definition}

\newtheorem{example}[theorem]{Example}
\theoremstyle{remark}
\newtheorem{remark}[theorem]{Remark}

\newtheorem{notation}[theorem]{Notation}
\newtheorem{notationassumptions}[theorem]{Notation and Assumptions}
\newtheorem{assumption}[theorem]{Assumption}
\newcommand{\Spec}{\mbox{Spec}}
\newcommand{\calA}{{\mathcal A}}
\newcommand{\calB}{{\mathcal B}}
\newcommand{\calC}{{\mathcal C}}

\newcommand{\calF}{{\mathcal F}}

\newcommand{\calG}{{\mathcal G}}

\newcommand{\calI}{{\mathcal I}}
\newcommand{\calJ}{{\mathcal J}}

\newcommand{\calO}{{\mathcal O}}

\newcommand{\calZ}{{\mathcal Z}}
\newcommand{\sZ}{\calZ}

\newcommand{\Q}{{\mathbb Q}}

\newcommand{\Z}{{\mathbb Z}}
\newcommand{\A}{{\mathbb A}}

\newcommand{\PP}{{\mathbb P}}

\newcommand{\pp}{{\mathfrak p}}
\newcommand{\qq}{{\mathfrak q}}

\newcommand{\Pp}{{\mathfrak P}}

\newcommand{\Ii}{{\mathfrak I}}
\newcommand{\Jj}{{\mathfrak J}}

\newcommand{\ord}{\mbox{ord}}

\usepackage{textcomp}

\usepackage{calrsfs}
\begin{document}
 \title{Turing degrees of isomorphism types of geometric objects}
\author{Wesley Calvert}
\address{Department of Mathematics\\
Mail Code 4408\\
Southern Illinois University\\
1245 Lincoln Drive\\
Carbondale, Illinois 62901}
\email{wcalvert@siu.edu}
\author{Valentina Harizanov}
\address{Department of Mathematics\\
George Washington University\\
Washington, DC 20052}
\email{harizanv@gwu.edu}
\author{Alexandra Shlapentokh}
\address{Department of Mathematics\\
East Carolina University\\
Greenville, NC 27858}
\email{shlapentokha@ecu.edu}
\date{\today}
\begin{abstract}

We initiate the computability-theoretic study of ringed spaces and
schemes.  In particular, we show that any Turing degree may occur as
the least degree of an isomorphic copy of a structure of these kinds.
We also show that these structures may fail to have a least degree.

\end{abstract}

\thanks{Calvert acknowledges partial support of the NSF grant
  DMS-0554841 and a Fulbright-Nehru Senior Research Fellowship, Calvert and Harizanov of the NSF grant DMS-110112, Harizanov of the NSF grant DMS-0904101, and Shlapentokh of
the NSF grant DMS-0650927, a grant from Templeton Foundation and an ECU Faculty Senate Summer Grant.}
\maketitle
\section{Introduction}

Hartshorne states the problem of classifying all algebraic varieties
up to isomorphism as one of the central guiding problems of algebraic
geometry, ``which are so difficult that one never expects to solve
them completely, yet which provide stimulus for a great amount of
work, and which serve as yardsticks for measuring progress'' \cite{Har}.

The present paper does not address this problem directly, but does
bear upon it.  Toward the problem
of classifying varieties up to isomorphism, it would be interesting to
know that we can associate each Turing degree with an isomorphism type
of varieties.  In this way, one could say that classifying varieties
up to isomorphism entailed a structure as complicated as the upper
semilattice of Turing degrees.

Interestingly, this encoding does not seem possible.  Once the field
is fixed, it is difficult to vary any undecidability properties among
varieties.  In the present paper, then, we work with structures
slightly richer than varieties: ringed spaces representing unions of
varieties, and schemes.

Among these structures, we do show that we can associate each Turing
degree to an isomorphism type.  This is interesting for two reasons:
\begin{enumerate}
\item It shows that the problem of classifying these structures up to
  isomorphism would produce a system of invariants at least as
  complicated as the upper semilattice of Turing degrees.
\item The particular association of degrees with structures happens in
  the undecidability of the basic algebraic and geometric features of
  the structures.  Thus, our results show that arbitrary levels of
  algorithmic ``information'' may be encoded in these structures.
\end{enumerate}

Notice that we are trying to demonstrate the great
diversity possible in these classes.  An important consequence of this
is that the structures we construct will, of necessity, be far outside
the range of examples normally considered.  To show that pathology is
possible in a category, one must construct pathological examples, and
all of the examples in this paper are, in some sense, pathological.

A broader goal of this paper is to initiate work in model theory and
computable model theory on these classes of objects, especially
schemes.  To date, while there has been work on varieties, compact
complex manifolds, and K\"ahler manifolds \cite{moosa,moosak,pillay}, we have been able to find
no work on schemes in model theory or computable model theory.

We believe that this does not owe to a lack of interest in
schemes, but rather to their reputation for difficulty.  Since the
calculation of degrees of isomorphism types is a central question of
computable model theory, it is our hope that our success with this
problem will encourage more work in this area.  In any case, Section
\ref{seclanguage} provides the logical infrastructure that any such work would
require: a determination of a language whose morphisms correspond
exactly to the classical homomorphisms.

More precisely, we study Turing degrees of isomorphism types
of structures from some well-known classes. This is a natural way,
introduced by Jockusch and Richter (see \cite{R}), of expressing the
algorithmic complexity of the structure. We consider only countable
structures for computable languages. The universe $A$ of an
infinite countable structure $\mathcal{A}$ can be identified with the set $%
\omega $ of all natural numbers. Furthermore, we often use the same symbol
for the structure and its universe. (For the definition of a language and a
structure see p.\ 8 of \cite{Marker}, and for a definition of a computable
language see p.\ 509 of \cite{Millar}.)

The two motivating constructions of the paper are of a union of varieties
(Section \ref{secunionvar}) and schemes (Section \ref{secschemes}).
In the remainder of the present section, we will describe some necessary algorithmic and geometric
background.  Section \ref{seclanguage} will describe a language for
treating varieties, schemes, and other ringed spaces as first order structures (in
the sense of model theory).  The remaining three sections (excluding
the appendix) will each develop an example of a class of structures
that admits (in a sense to be made precise in Section
\ref{ssecturing}) the encoding of arbitrary Turing degrees and of a
minimal pair of Turing degrees.  Section \ref{secunionweak} explores
the class of unions of subspaces of a fixed variety, under a weak
topology.  Section \ref{secunionvar} constructs a union of varieties,
and section \ref{secschemes} constructs schemes.

\subsection{Definitions, and Statement of Main Results}\label{ssecturing}

We say that a set $X$ is \emph{Turing reducible to} (\emph{computable in}) a
set $Y$, in symbols $X\leq _{T}Y$, if $X$ can be computed by an algorithm
with $Y$ in its oracle. Turing reducibility is the more basic notion, in
terms of which Turing degree is defined. We say that the sets $X$ and $Y$
are \emph{Turing equivalent}, or have the same \emph{Turing degree}, if $%
X\leq _{T}Y$ and $Y\leq _{T}X$. We use $\equiv _{T}$ for Turing equivalence.
We also write $deg(X)=deg(Y)$ or $Y\in deg(X)$ instead of $X\equiv _{T}Y$.
(Detailed information about Turing degrees and their structure can be found
in \cite{Rogers} and \cite{S}.)

When measuring complexity of structures, we identify them with their atomic
diagrams. The atomic diagram of a structure $\mathcal{A}$ is the set of all
quantifier-free sentences in the language of $\mathcal{A}$ expanded by
adding a constant symbol for every $a\in A$, which are true in $\mathcal{A}$. The \emph{Turing degree} of $\mathcal{A}$, $deg(\mathcal{A})$, is the
Turing degree of the atomic diagram of $\mathcal{A}$. Hence, $\mathcal{A}$ is
\emph{computable } iff $deg(\mathcal{A})=\mathbf{0}$.
(Some authors call a structure computable if it is only isomorphic to a
computable one.) We also say that a set or a procedure is \emph{computable}
(effective) \emph{relative to} $\mathcal{B}$, or \emph{in} $\mathcal{B}$,
if it is computable relative to the atomic diagram of $\mathcal{B}$\textit{.}

The \emph{Turing degree spectrum} of a countable structure $\mathcal{A}$ is
\begin{equation*}
DgSp(\mathcal{A})=\{\deg (\mathcal{B}):\mathcal{B}\cong \mathcal{A}\}\text{.}
\end{equation*}%
A countable structure $\mathcal{A}$ is \emph{automorphically trivial} if
there is a finite subset $X$ of the domain $A$ such that every permutation
of $A$, whose restriction on $X$ is the identity, is an automorphism of $%
\mathcal{A}$. If a structure $\mathcal{A}$ is automorphically trivial, then
all isomorphic copies of $\mathcal{A}$ have the same Turing degree. It was
shown in \cite{HM} that if the language is finite, then that degree must be $%
\mathbf{0}$. On the other hand, Knight \cite{K} proved that for an
automorphically nontrivial structure $\mathcal{A}$, we have that $DgSp(%
\mathcal{A})$ is closed upwards, that is, if $\mathbf{b}\in DgSp(\mathcal{A}%
) $ and $\mathbf{d}>\mathbf{b}$, then $\mathbf{d}\in DgSp(\mathcal{A})$.
Hirschfeldt, Khoussainov, Shore, and Slinko \cite{HKSS} established that for
every automorphically nontrivial structure $\mathcal{A}$, there is a
symmetric irreflexive graph, a partial order, a lattice, a ring, an integral
domain of arbitrary characteristic, a commutative semigroup, or a $2$-step
nilpotent group whose degree spectrum coincides with $DgSp(\mathcal{A})$.

Since the Turing degree of a structure is not invariant under isomorphisms,
Jockusch and Richter introduced the following complexity measures of the
isomorphism type of a structure.

\begin{definition}
$($i$)$ The \emph{Turing degree of the isomorphism type} of $\mathcal{A}$,
if it exists, is the least Turing degree in $DgSp(\mathcal{A})$.

$($ii$)$ Let $\alpha $ be a computable ordinal. The $\alpha $\emph{th jump
degree} of a structure $\mathcal{A}$ is, if it exists, the least Turing
degree in%
\begin{equation*}
\{deg(\mathcal{B})^{(\alpha )}:\mathcal{B}\cong \mathcal{A}\}\text{.}
\end{equation*}
\end{definition}

\noindent Obviously, the notion of the $0$th jump degree of $\mathcal{A}$
coincides with the notion of the degree of the isomorphism type of $\mathcal{%
A}$. (A general discussion of the jump operator can be found in 13.1 of \cite%
{Rogers}, and Chapter III\ of \cite{S}.)

In \cite{R} Richter proved that if $\mathcal{A}$ is a structure without a
computable copy and satisfies the effective extendability condition
explained below, then the isomorphism type of $\mathcal{A}$ has no degree.
Richter's result uses a minimal pair construction. Distinct nonzero Turing
degrees $\mathbf{a}$ and $\mathbf{b}$ form a \emph{minimal pair} if
\begin{equation*}
(\mathbf{c}\leq \mathbf{a} \wedge \mathbf{c}\leq \mathbf{b)}\Rightarrow
\mathbf{c}=\mathbf{0}\text{.}
\end{equation*}%
(See \cite{S} for the minimal pair method.) A structure $\mathcal{A}$
satisfies the effective extendability condition if for every finite
structure $\mathcal{M}$ isomorphic to a substructure of $\mathcal{A}$, and
every embedding $\sigma $ of $\mathcal{M}$ into $\mathcal{A}$, there is an
algorithm that determines whether a given finite structure $\mathcal{N}$
extending $\mathcal{M}$ can be embedded into $\mathcal{A}$ by an embedding
extending $\sigma $. Richter \cite{R} also showed that every linear order,
and every tree as a partially ordered set, satisfy the effective
extendability condition.  More recent results are described in the
introduction to \cite{CHS}, where using results of  Richter \cite{R},   the authors proved the following:

\begin{theorem}
\label{thm:2}
If $\mathcal{C}$  any of the classes listed below, then for any Turing degree $\mathbf{d}$ there is a member  of $\mathcal{C}$ whose isomorphism type has Turing degree
  $\mathbf{d}$.  There is also a member of $\mathcal{C}$ whose   isomorphism type has no Turing degree.
\begin{enumerate}
\item The class of algebraic extensions of any fixed computable finitely  generated field
\item The family of sub-rings $R$ of any fixed computable finitely   generated field $K$, such that the fraction field of $R$ is $K$
\item The class of torsion-free Abelian groups of any fixed finite  rank
\end{enumerate}
\end{theorem}

The main results of the present paper add structures of a very
different kind to this list.

\begin{theorem}[Main Result] If $\mathcal{C}$ is any of the classes listed
  below, then for any Turing degree $\mathbf{d}$ there is a member  of $\mathcal{C}$ whose isomorphism type has Turing degree
  $\mathbf{d}$.  There is also a member of $\mathcal{C}$ whose
  isomorphism type has no Turing degree.
\begin{enumerate}
\item Ringed spaces corresponding to unions of subvarieties of certain
  fixed varieties
\item Ringed spaces corresponding to unions of varieties
\item Schemes over a fixed field.
\end{enumerate}
\end{theorem}

The results in \cite{CHS} were proved using the two following
cumbersome, but not technically difficult frameworks.

\begin{theorem}[\cite{CHS}]%
\label{thm:int1}%
 Let $\mathcal{C}$ be a class of countable structures in a
computable language $L$, closed under isomorphisms. Assume that there is a
computable sequence $\{\mathcal{A}_{i},i\in \omega \}$ of computable ( possibly infinite ) structures in $\mathcal{C}$ satisfying the following
conditions.

\begin{itemize}
\item There exists a finitely generated structure $\mathcal{A}\in \mathcal{C}
$ such that for all $i\in \omega $, we have that $\mathcal{A}\subset
\mathcal{A}_{i}$.

\item For any $X\subseteq \omega $, there is a structure $\mathcal{A}_{X}$
in $\mathcal{C}$ such that $\mathcal{A}\subset \mathcal{A}_{X}$ and
\begin{equation}
\label{eq:int1}
\mathcal{A}_{X}\leq _{T}X\text{,}
\end{equation}%
and for every $i\in \omega $, there exists an embedding $\sigma $ such that
\begin{equation*}
\sigma :\mathcal{A}_{i}\hookrightarrow \mathcal{A}_{X}\text{, \ }\sigma _{|%
\mathcal{A}}=id\text{,}
\end{equation*}%
iff $i\in X$.

\item Suppose that a structure $\mathcal{A}_{X}$ is isomorphic to some structure $\mathcal{B}$ under isomorphism
\[
\tau :\mathcal{A}_{X}\longleftrightarrow \mathcal{B},
\]
and let a pair of structures $\mathcal{A}_{i},\mathcal{A}_{j}$ be such that exactly one of them embeds in $\mathcal{B}$ via $\sigma $ with $(\tau ^{-1}\circ \sigma )_{|\mathcal{A}}=id$. In this case there is a uniformly effective procedure with oracle $\mathcal{B}$ for deciding which of the two structures embeds in $\mathcal{B}$.
\end{itemize}

Under these assumptions,  for every Turing degree $\mathbf{d}$, there is a structure in $\mathcal{C}$ whose isomorphism type has degree $\mathbf{d}$.
\end{theorem}

The second framework theorem has hypotheses that look similar to those of the
first, but has a different conclusion.  For the new hypotheses, we
need the following standard definition \cite{Rogers}.

\begin{definition} We say that $X \leq_e Y$ if we can uniformly pass
  from any enumeration of $Y$ to an enumeration of
  $X$.\end{definition}

By replacing Turing reducibility with enumeration reducibility in two
key points, we obtain the following statement:

\begin{theorem}[\cite{CHS}]
\label{thm:int2}
Let $\mathcal{C}$ be a class of countable structures in a computable language $L$%
, closed under isomorphisms. Assume that there is a computable sequence $\{%
\mathcal{A}_{i},i\in \omega \}$ of computable $($possibly infinite$)$
structures in $\mathcal{C}$ satisfying the following conditions.

\begin{itemize}
\item There exists a finitely generated structure $\mathcal{A}\in \mathcal{C}
$ such that for all $i\in \omega $, we have that $\mathcal{A}\subset
\mathcal{A}_{i}$.

\item For any $X\subseteq \omega $, there is a structure $\mathcal{A}_{X}$
in $\mathcal{C}$ such that $\mathcal{A}\subset \mathcal{A}_{X}$ and
\begin{equation}
\label{eq:int2}
\mathcal{A}_{X}\leq _{e}X\text{,}
\end{equation}%
and for every $i\in \omega $, there exists an embedding $\sigma $ such that
\begin{equation*}
\sigma :\mathcal{A}_{i}\hookrightarrow \mathcal{A}_{X}\text{, \ }\sigma _{|%
\mathcal{A}}=id\text{,}
\end{equation*}%
iff $i\in X$.

\item If a structure  $\mathcal{A}_{X}$ is isomorphic to some structure $\mathcal{B}$ under isomorphism
\[
\tau :\mathcal{A}_{X}\longleftrightarrow \mathcal{B},
\]
then from any enumeration of $\mathcal{B}$, we can effectively enumerate those $i$ for which $\mathcal{A}_{i}\hookrightarrow \mathcal{B}$
under an embedding $\sigma $ with $(\tau ^{-1}\circ \sigma )_{|\mathcal{A}}=id$.
\end{itemize}
Under these assumptions, there is a structure $\mathcal{A}_{X}$ in $\mathcal{C}$ whose isomorphism type has no Turing degree.
\end{theorem}

These same two results will give the computability-theoretic basis for
the results of the present paper.  The major challenge addressed in
the remaining pages is showing that the hypotheses of these theorems
are satisfied by the classes listed in the Main Result.

\subsection{Geometric Preliminaries}
We start with defining the main objects of the paper and we will proceed from the more general (pre-sheaves) to more specific (schemes).
Most of the following definitions can be found in \cite{Har}.
\begin{definition}[A pre-sheaf and a sheaf of rings]
\label{def:sheaf}
 Let $X$ be a topological space.  A  pre-sheaf ${\mathcal F}$ of rings on $X$ consists of the data
 \begin{enumerate}%
\item for every open subset $U \subseteq X$, a ring ${\calF}(U)$, and \item for every inclusion $V \subseteq U$ of open subsets of $X$ a morphism of rings $$\rho_{U,V} : {\calF}(U) \rightarrow {\calF}(V),$$
\end{enumerate}%
subject to the following conditions:
\begin{enumerate}%
 \item ${\calF}(\emptyset) = \{0\}$, %
\item $\rho_{U,U}$ is the identity map ${\calF}(U) \longrightarrow
  {\calF}(U)$, and
\item If $W \subseteq V \subseteq U$ are three
open sets, then $\rho_{U,W} = \rho_{V,W}\circ\rho_{U,V}$.
\end{enumerate}

A pre-sheaf ${\mathcal F}$ on a topological space $X$ is a {\it sheaf} if it satisfies the following supplementary
conditions:
\begin{enumerate}
\item[(4)]   If $U$ is an open set, $\{V_i\}$ is an open covering of $U$, and  $s \in {\calF}(U)$ is an element such that $\rho_{U,V_i}(s)=0$ for
all $i$, then $s=0$.
\item[(5)]  If $U$ is an open set,  $\{V_i\}$ is an open covering of  $U$, and for each $i$ there exists an  element $s_i \in {\calF}(V_i)$    with the property that for any  $j$, we have
\[
\rho_{V_i,V_i\cap V_j}(s_i)=  \rho_{V_j,V_i\cap V_j}(s_j),
\]
then there exists $s \in {\calF}(U)$ such that for each $i$ it is the case that  $\rho_{U,V_i}(s)= s_i$.
\end{enumerate}%
\end{definition}

\begin{definition}[Direct image of a sheaf]
Let $f: X \longrightarrow Y$ be a continuous map of topological spaces. For any sheaf ${\calF}$ on $X$, we define the direct image
sheaf $f_*{\calF}$ on $Y$ by setting $(f_*{\calF})(V)={\calF}(f^{-1}(V))$ for any open subset $V$ of $Y$.
\end{definition}

\begin{definition}[Inverse image of a sheaf]
Let $f: X \longrightarrow Y$ be a continuous map of topological spaces. For any sheaf ${\calG}$ on $Y$, we define the inverse  image
sheaf $f^{-1}{\calG}$ on $X$ by setting
\[
(f^{-1}{\calG})(U)= \lim_{f(U)\subseteq V \subseteq Y, V \mbox{\tiny open}}\calG(V)
\]
for any open subset $U$ of $X$.
\end{definition}

\begin{definition}[Morphism of sheaves]
If ${\calF}$ and ${\calG}$ are sheaves on a topological space $X$, then  a morphism  $\phi : {\calF} \rightarrow {\calG}$ consists of the morphism of rings $$\phi(U): {\calF}(U) \rightarrow {\calG}(U)$$ for each open set $U$ of $X$, such that whenever $V \subseteq U$  holds the diagram commutes.
\begin{center}
$\begin{CD}
{\mathcal{F}(U)}  @>{\phi(U)}>>  {\mathcal{G}(U)}\\
@V{\rho_{U,V}}VV    @VV{\rho'_{U,V}}V\\
{\mathcal{F}(V)}  @>>{\phi(V)}>  {\mathcal{G}(V)}\\
\end{CD}$
\end{center}
\end{definition}

\begin{definition}[Restriction sheaf]
Let $Z$ be a subset of a topological space $V$ under restriction topology.  If $\calF$ is a sheaf on $V$ and  $i : Z \longrightarrow V$ is the inclusion map, then the inverse image of $\calF$ under $i$, that is $i^{-1}\calF$ is called the restriction of $\calF$ to $Z$.

\end{definition}
Below  are definitions of  ringed spaces, varieties and schemes.
\begin{definition}[Ringed space]%
A ringed space is a pair $(X, {\calO}_X)$ consisting of a topological space  $X$ and a sheaf of rings ${\calO}_X$ on $X$.
An morphism of ringed spaces from $(X,{\calO}_X)$ to $(Y, {\calO}_Y)$ is a pair

$(f,f^{\#})$ of a continuous map $f: X \longrightarrow Y$ and an morphism
$f^{\#}: {\calO}_Y \longrightarrow f_*{\calO}_X$ of sheaves of rings on $Y$.
\end{definition}%
 \begin{definition}[Stalk]%
If  ${\calF}$ is a sheaf on a topological space $Y$,  $Q \in Y$ and  $V$ ranges over all
open neighborhoods of $Q$, then  the direct limit  $\lim\limits_{\overrightarrow{V}}{\calF}(V) $ is called
the stalk ${\calF}_Q$ of ${\calF}$ .
\end{definition}%
\begin{definition}[Locally ringed spaces]%
A ringed space $(X, {\calO}_X)$ is a locally ringed space if for each point $P \in X$, the stalk ${\calO}_{X,P}$ is a local ring, i.e. a ring with a unique maximal ideal.
A morphism of locally ringed spaces is a morphism $(f,f^{\#})$ of ringed spaces, such that for each point $P \in X$, the induced
map of local rings $f^{\#}_P : {\calO}_{Y,f(P)} \longrightarrow \calO_{X,P}$ is a homomorphism of local rings.
\end{definition}%
\begin{definition}[Spectrum of a ring as a ringed space]
Let $R$ be a ring and let $\Spec R$ be the set of all of its prime ideals.  A subset $\calJ$ of $\Spec R$ is closed under a Zariski topology if there exists a non zero ideal $I$ contained in every prime ideal of $\calJ$.  Given an open subset $O \subseteq \Spec R$ we let $\calF(O) = \bigcap\limits_{\pp \in O}R_{\pp}$, where $R_\pp$ is the localization of $R$ at $ \pp$ and the intersection is taken inside $\coprod\limits_{\pp \in \Spec R}R_{\pp}$.
\end{definition}

\begin{definition}[Scheme]%
An affine scheme is a locally ringed space $(X, {\calO}_X)$ which is isomorphic (as a locally ringed space) to the spectrum of
some ring. A scheme is a locally ringed space $(X, {\calO}_X)$ in which every point has an open neighborhood $U$ such that
topological space $U$, together with restricted sheaf ${\calO}_{X|U}$ is an affine scheme. We call $X$ the underlying topological
space of the scheme $(X,{\calO}_X)$, and ${\calO}_X$ its structure sheaf.  A morphism of schemes is a morphism of schemes as locally ringed spaces.
\end{definition}%

\begin{definition}[Scheme over a scheme]%
\label{def:over}
Let $S$ be a fixed scheme.  A scheme over $S$ is a scheme $X$ together with  a morphism $X
\rightarrow S$.  If $X, Y$ are schemes over $S$, a morphism from $X$ to $Y$ as schemes over $S$ is a
morphism $X \longrightarrow Y$ compatible with the given morphism to $S$.
\end{definition}%

We will also use a similar construction for ringed spaces: a ringed space over a ringed space defined in an analogous manner.

We finally make several  observations which will be helpful in the sequel.
The following proposition can be found in Section 15.10.1 of  \cite{Eisenbud}.
\begin{proposition}
\label{prop:idealcomput}
Let $k$ be a computable field.  If $I$ is an ideal of $R=k[x_1,\ldots, x_n]$ generated by $f_1,\ldots,f_m \in R$, then $I$ is a computable subset of $R$ and consequently the coordinate ring of any variety over $k$ is computable.
\end{proposition}
Then next lemma is also a well-known fact which we state here for the convenience of the reader.
\begin{lemma}
\label{le:closedset}
Let $k$ be a computable field, $V$ a  variety over $k$
  given by a system of polynomial equations, and let $f \in k(V)$;
  that is, in the function field of $V$.  Then given an open $U \subseteq
  V$, we may effectively decide whether $f \in k(U)$, i.e. whether $f$ is defined on $U$.
  \end{lemma}

\begin{proof}
Assume $V \subset k^n$. If $f =\frac{f_1}{f_2}$, where $f_1, f_2 \in
k[x_1,\ldots,x_n]$  are relatively prime, and the equivalences classes
of $f_1, f_2$ are denoted by $[f_1], [f_2] \in k[V]$, then $[f] \in
k(U)$ if and only if the set of zeros of $f_2$, denoted by $V(f_2)$,
is contained in $U^c$, the complement of $U$ in $V$ .  Deciding
whether $V(f_2)$ is contained in $U^c$ is equivalent to deciding
whether $I(U^c)$, the ideal of polynomials classes in $k[V]$ having a zero at every point of $U^c$, is contained in the ideal generated by the equivalence class of $[f_2]$.  Let $\{[v_i]\}$ be a (finite) set of generators for $I(U^c)$. Divide each $v_i$ of them by $f_2$, and check if the remainder is zero in $k[V]$.

\end{proof}
In a similar fashion one can prove the following lemma.
\begin{lemma}
\label{le:closedset1}
Let $k$ be a computable field, $V$ a  variety over $k$
  given by a system of polynomial equations, and let $f \in k(V)$;
  that is, in the function field of $V$.  Then given a system of equations  defining a closed set $W \subseteq
  V$, we may effectively decide whether $f \equiv 0$ on $W$.
  \end{lemma}

We will also use the following proposition.
\begin{lemma}
\label{le:infnum}
Let $V$ be a variety over a countable field $k$, and let $\{W_i, i \in \Z_{>0}\}$ be a countably infinite set of pairwise disjoint subvarieties.  If $f \in k(V)$, then $f$ is not identically zero on $\bigcup_{i \in \Z_{>0}}W_i$ unless $f \equiv 0$ on $V$.
\end{lemma}
\begin{proof}
Write $f=\frac{f_1}{f_2}$, where $f_1, f_2 \in k[V]$, which is a Noetherian ring.  Since each $W_i$ is a subvariety, the ideal $I_i$ of all functions zero on $W_i$ is a prime ideal of $k[V]$ and for any $i \not =j$ we have that $I_i \not =I_j$.  If $V(f_1) \supset W_i$, then $I(f_1) \subset I_i$.  Since a ring is Noetherian, any non-zero ideal can be contained in only finitely many prime ideals.  Thus, $f$ cannot be identically zero on every $W_i$ without being identically zero on $V$.
\end{proof}
Finally, we will need the following fact.
\begin{lemma}
Let $V$ be a variety over a computable algebraically closed field $k$ embedded in an $n$-dimensional affine pace. Let $f_1, f_2 \in k[x_1,\ldots,x_n]$ be two polynomials such that $f_2 \not \in I(V)$.  In this case there is a recursive procedure uniform in equations of $V$ to determine if $\frac{f_1}{f_2}$ is a constant function in $k(V)$.
\end{lemma}
\begin{proof}
First of all we note that whether or not $f_2 \in I(V)$ can be established algorithmically by Proposition \ref{prop:idealcomput}. Assuming  $f_2 \not \in I(V)$, Let $ \bar c \in V$ be such that $f_2(\bar c) \not =0$.   (Such $\bar c$ can be located by a search of $V$, which is computable).  Now observe that $\frac{f_1}{f_2}=\mbox{const}$ on $V$ if and only $f_1(\bar x)f_2(\bar c)-f_2(\bar x)f_1(\bar c) \in I(V)$
\end{proof}
\section{A Language for Geometric Structures}\label{seclanguage}
\setcounter{equation}{0}
A major problem preliminary to the present work is finding a way to
talk about geometric objects.  We are not the first to introduce a first order model-theoretic
language for these structures.  Zilber, Moosa
\cite{moosa}, and Pillay \cite{pillay}, for instance, use the
following language for compact complex manifolds, and note that
varieties are a special case.

\begin{definition} Let $L_1$ be the language with countably many predicates in
  each arity.  Let $V$ be a variety.  We interpret the predicates to
  represent each Zariski closed set in $V$. \end{definition}

We  follow along the same path, but since we consider arbitrary ringed spaces we will have to add means to describe the rings corresponding to the open sets.  

In the following definition, the reader may wish to refer back to the
definition of a pre-sheaf in Section 1.2.  The ringed spaces and schemes described below will have countable
underlying topological spaces, with countably many open and closed
sets and countably many elements in each ring corresponding to an open
subset of the topological space.  Thus we propose to use a language
$L_R$ which has the following components:

\begin{itemize}
\item For each open set a predicate for membership in the set.
\item A ring language
for each ring.
\item For each open set $U$ a predicate $F_U(x)$ with the interpretation that $F_U(x)$ is true if and only if $x \in \calF(U)$.
\item For each pair of open sets $U \subset V$ a function $R_{U, V}(x,y)$, with the interpretation that $R_{U, V}(x,y)$ is true if and only if $x \in \calF(U), y\in \calF(V)$ and $\rho_{U,V}(x)=y$.
\end{itemize}

\begin{proposition}
 Let $V$ and $U$ be ringed spaces.
 Then the following are identical:
\begin{enumerate}
\item The $L_R$-morphisms from $V$ to $U$,
\item The morphisms from $V$ to $U$ in the category of ringed spaces.
\end{enumerate}
\end{proposition}

\begin{proof}
Let $\varphi: V \to U$ be an $L_R$ morphism.  Since $\varphi$ respects
open sets (in both directions), it must be continuous on the underlying
topological space.  Since $\varphi$ respects both membership in
$\calF(U)$ for each $U$ and the ring structure of that set, as well as
respecting inclusions via the functions $R_{U,V}$, we know that
  $\varphi$ must be a morphism of sheaves.  The other direction is symmetric.
\end{proof}

\section{Unions of ``Subspaces''}\label{secunionweak}
In this section we will describe the first of three general classes of ringed spaces which satisfy the conditions of Theorems \ref{thm:int1} and \ref{thm:int2}.  Below we construct ringed spaces made out of unions of closed irreducible subsets.  We start with a general fact about topological spaces.
\begin{lemma}
\label{le:intersec}
If $Z$ is a topological space, then it is  irreducible (i.e. not equal
to a union of two proper closed subsets) if and only if  for any
non-empty open sets $U$ and $T$, the intersection $U\cap T$ is not empty.
\end{lemma}
\begin{proof}
Let $U, T$ be as in the statement of the lemma and without loss of
generality assume that neither is equal to the whole space.  Let $U^C,
T^C$ be the complements of $U$ and $T$, respectively.  Since $Z$ is
irreducible, given our assumptions, we have that $U^C \cup T^C \not =
V$.  Thus, $U\cap T \not = \emptyset$.  Conversely, if there exist
closed subsets $A,B \subseteq Z$ such that neither is the whole space and $Z=A \cup B$, then $A^C \cap B^C= \emptyset$ while neither $A^C \not = \emptyset$ nor $B^C \not = \emptyset$.
\end{proof}
We now define a class of ringed spaces of particular interest to us.
\begin{definition}[The $\calZ$-class]
A ringed space $(Z, O_Z)$ will belong to the $\calZ$-class if the following conditions are satisfied.
\begin{enumerate}
\item The topological space $Z$ is countable.
\item The number of open/closed subsets of $Z$ is countable.
\item $\calF(Z)=R(Z)$ is a countable ring of functions from $Z$ to some field $k$.  The ring $R(Z)$ is an integral domain and is also a $k$-algebra.  The fraction field of $R(Z)$ is denoted by $K(Z)$.
 \item For any open $U \subseteq Z$ we have that  $\calF(U) \subset
   K(Z)$ and $\calF(U)$ contains functions of $K(Z)$ defined on $U$ (but not necessarily all the functions defined on $U$).  A function $f \in K(Z)$ is defined on $U$ if $f = \frac{f_1}{f_2}$, with $f_1, f_2 \in R(Z)$ and $f_2(z) \not = 0$ for any $z \in U$.
\item  For  $U \subseteq V$  open sets, we have that   $\calF(V)\subseteq \calF(U)$ and $\rho_{U,V}$ is defined to be the inclusion map.
\item For any collection $A$ of open subsets of $Z$ we have that $\bigcap_{U \in A}\calF(U)=\calF(\bigcup_{U \in A}U)$.
 \end{enumerate}
\end{definition}

\begin{definition}[The zero-set condition]
Let $A \subseteq Z$ be closed and irreducible and satisfy the following condition.
\begin{center}
{\it If $U$ is an open subset of $Z$, $f \in \calF(U)$ and $f$ is identically zero on $A \cap U$, then $f$ is defined for every point of $A$ and its value is equal to zero at every point of $A$.}
\end{center}
In this case we will say that $A$ satisfies the {\it zero-set condition}.
\end{definition}
We will also need the following objects:
\begin{notation}
Let $G(Z)$ be the smallest ring inside $K(Z)$ containing $\calF(U)$ for any open set $U$ of $Z$.  (Alternatively, this is the ring generated by elements of all $\calF(U)$ when considered as elements of $K(Z)$.)
\end{notation}
We now construct a new type of a ringed space out of subsets of $Z$.
\begin{definition}[A Union ringed subspace]
 For each positive integer $i$, let $Z_i \subseteq Z$ be closed,  irreducible, and satisfy the zero-set condition  as a subset of $Z$. Let $Y = \bigcup Z_i$ and consider a topology on $Y$ where open sets are of the form $U_Y=U \cap Y,$ with $U$ being an open subset of $Z$  such that for all  $i$ we have that $U_i= Z_i \cap U$ is a non-empty (open) subset of $Z_i$ (in the relative topology inherited from $Z$).  Let
 \[
 \calI(Y) = \{f \in G(Z): f|_{Y} \equiv 0\},
 \]
 let $\calF(U_Y)$ be the ring generated by $\bigcup_{U \supseteq U_Y,U \mbox{ open }}\calF(U)$ and assume that $\calI(Y) \subset \calF(Y)$.   Now define $\calF_Y(U_Y) = \calF(U_Y)/\calI(Y)$.
 \end{definition}
 \begin{remark}
 \label{rem:I}
Using the notation above, we make the following observations.
\begin{enumerate}
\item \label{it:1} {\it If $U_Y$ is an open subset of $Y$ and $f \in \calF(U_Y)$, then $f \in \calF(W)$ for some open $Z$-subset $W$ containing $U_Y$.}  Indeed, by definition, $\calF(U_Y)$ is    generated by elements of all rings $\calF(U)$ such that $U$ is an open set containing $U_Y$.  Thus $f = \sum_{i_1,\ldots,i_k}f_{i_1}\ldots f_{i_k}$, where each $f_{i_j} \in \calF(U_{i_j})$, with $U_{i_j}$ an open set containing $U_Y$.  If $$W=\bigcap_{i_1,\ldots, i_k} U_{i_j},$$ then $W$ is an open subset of $Z$ containing $U_Y$ with $W \subseteq U_{i_j}$ for each $i_j$, and by a property of the $\calZ$-class, we have that each $\calF(U_{i_j}) \subseteq \calF(W)$ implying each $f_{i_j} \in \calF(W)$.  Thus $f \in \calF(W)$.

\item \label{it:2} {\it If $W_Y \subseteq U_Y$ are open subsets of $Y$, then $\calF(U_Y)\subseteq \calF(W_Y)$.}  This follows from the fact any open subset $U$ of $Z$ containing $U_Y$ will contain $W_Y$.

\item  {\it $\calI(Y) = \{f \in \calF(U_Y): f|_{U_Y} \equiv 0 \}$. } Observe that clearly
\[
\calI \subseteq \{f \in \calF(U_Y): f|_{U_Y} \equiv 0 \},
\]
and at the same time if $f \in \calF(U_Y)\in G(Z)$ is identically zero on $U_Y$, it is identically zero on an open subset (in the relative topology inherited from $Z$)  of each $Z_i$ and, since $Z_i$ satisfy the zero-set condition,  $f$ is therefore  identically zero on each $Z_i$ and thus on $Y$, implying it is in $ \calI(Y)$.
\end{enumerate}
 \end{remark}
  \begin{proposition}
  \label{prop:well-def}
$ (Y, O_Y)$ where $O_Y =\{\calF_Y(U_Y), U_Y \mbox{ an open subset of } Y\}$ is a well-defined ringed space.
  \end{proposition}
  \begin{proof}
 First we show that the topology of $Y$ is well defined. Observe that under the relative topology inherited from $Z$, by Lemma \ref{le:intersec}, no $Z_i$ has two non-empty open subsets with empty intersection.  Thus, if $U, V$ are two open subsets of $Z$ such that $U\cap Z_i \not = \emptyset$ and $V\cap Z_i \not = \emptyset$, we also have that $U\cap V\cap Z_i \not =\emptyset$.     Thus, $U\cap V$ will have a non-empty intersection with each $Z_i$.  Therefore the proposed class of open subsets of $Y$ is closed under finite intersections  while it is obviously closed under arbitrary unions. This shows that our topology is well-defined.

 Observe that we can interpret the equivalence classes of $\calF_Y(U_Y)$ as functions from $U_Y$ to $k$.  We will denote a class of a function $f \in \calF(U_Y)$ by $[f]$.  By Remark \ref{rem:I}, Part \ref{it:2} for any  $U_Y \subseteq V_Y$ open subsets of   $Y$ we have $\calF(V_Y) \subset \calF(U_Y)$ and therefore
 \[
 \calF_Y(V_Y)\subseteq \calF_Y(U_Y) \subset G(Z)/\calI(Y).
 \]
 Thus we can define a map $\rho_{U_Y,V_Y}: \calF_Y(U_Y)\longrightarrow \calF_Y(V_Y)$ to be the inclusion map.  Consequently, the collection of rings $\{\calF_Y(U_Y): U_Y \subset Y, U_Y \mbox{ open}\}$ produces a pre-sheaf structure on $Z$.

Suppose now that $U_Y= \bigcup U_{Y,j}$ is an open covering of an open set in the topology of $Y$.  If for some $j$, for some $S \in \calF_Y(U_Y)$  we have that $\rho_{U_Y, U_{Y,j}}(S) = 0$, then $S=0$ since $\rho_{U_Y, U_{Y,j}}$ is the inclusion map.

 Now assume we have two classes of functions $S_j \in \calF_Y(U_{Y,j})$ and $S_m \in \calF_Y(U_{Y,m})$,  and
 \[
 \rho_{U_{Y,j},U_{Y,j}\cap U_{Y,m}}(S_j)=\rho_{U_{Y,m},U_{Y,j}\cap U_{Y,m}}(S_m).
\]
Since $ \rho_{U_{Y,j},U_{Y,j}\cap U_{Y,m}}$ and  $\rho_{U_{Y,m},U_{Y,j}\cap U_{Y,m}}$ are inclusion maps, we must conclude that $S_j =S_m$ as functions in $G(Z)/\calI(Y)$.  Therefore $S_j \in \calF_Y(U_{Y,j})\cap \calF_Y(U_{Y,m})$.   If we fix $j$ and let $m$ vary, we will conclude that $S_j$ is in fact in  $\bigcap_{m \in \Z_{>0}}\calF_Y(U_{Y,m})$. Now let $f_j \in S_j$ be an element of $G(Z)$.  In this case, we have that
\[
f_j \in  \bigcap_{m \in \Z_{>0}}\calF(U_{Y,m}).
\]
Further, from Remark \ref{rem:I}, Part \ref{it:1} for some sequence of open subsets $U_m$ of $Z$ with each  $U_m$ satisfying $U_{Y,m}=U_m \cap Y$ we have that
\[
f_j \in  \bigcap_{m \in \Z_{>0}}\calF(U_{m})= \calF(\bigcup_m U_m),
\]
where the last equality followed by the assumption that $Z$ is a $\calZ$-class ringed space.  Now since $\bigcup U_{Y,m} = U_Y$ we also have that $U_Y \subseteq \left (\bigcup_m U_m \right ) \cap Y =W_Y$ -- an open subset of $Y$.  Since $\calF(W_Y) \subseteq \calF(U_Y)$ by Remark \ref{rem:I}, Part \ref{it:2}, we conclude that $f_j \in \calF(U_Y)$ and thus $S_j \in \calF_Y(U_Y)$.
\end{proof}
\begin{remark}
Open sets $U$ intersecting each $Z_i$ always exist.  For example, $U$ can be the whole space.  Of course, if this is the only such set, the topology will be rather weak.
\end{remark}

\begin{remark}
\label{rem:justone}
 Let $Z_i$ be as above and consider $Y=Z_i$.   In this case the topology of $Y$ is simply the relative topology of $Z_i$ inherited from $Z$.  Thus $(Y=Z_i, O_Y=O_{Z_i})$ defined as above is a ringed space.
\end{remark}
 We now define morphisms between two unions of ringed spaces.  We first define a class of morphisms between ringed spaces in general and then address our particular situation.
 \begin{definition}[Composition Morphism]
 Let $S, T$ be ringed spaces of functions into some field $k$.  If $\phi : S \longrightarrow T$ is a morphism of ringed spaces, then $\phi$ will be called a {\it composition morphism} if for any open set $O \subset T$ we have that $\phi^{\#}: \calF(O) \longrightarrow \calF(\phi^{-1}(O))$ is defined by $f \mapsto f \circ \phi$.
 \end{definition}

\begin{definition}%
\label{def:subvarieties}
Let $Z, V $ be elements of the $\calZ$-class. Let $Z_i \subset Z, V_i \subset V$ be closed irreducible subsets, and finally let $Y = \bigcup_{i=1}^{\infty}  Z_i, W=\bigcup_{i=1}^{\infty}  V_i$, be two union ringed subspaces.  A map $\phi: Y \longrightarrow W$ will be called a morphism of union ringed subspaces  if $\phi$ is a restriction of a composition morphism $\bar \phi: Z  \longrightarrow V$  such that for any subset $S_i$ open in the relative topology  of  $V_i$ (inherited from $V$) we have that
\[
\phi^{-1}(S_i)\cap Y_j =\emptyset \iff \phi^{-1}(V_i) \cap Y_j=\emptyset.
\]
    For any open subset $O_W \subset W$ the ring homomorphism $$\phi^{\#}(O_W) : \calF_W(O_W) \longrightarrow \calF_Y(\phi^{-1}(O_W))$$ will be defined as follows.  A function $f \in \calF_W(O_W)$ will be mapped to a function $$f\circ \phi \in \calF_W(\phi^{-1}(O_W)).$$

An isomorphism $\phi$ between $Y$ and $W$ is an injective morphism whose inverse is also a morphism.  An embedding of $Y$ into $W$ is an injective morphism of $Y$ into $W$.
\end{definition}%
We now show the map $\phi$ as defined above is a morphism of sheaves of rings.  First we observe the following.
\begin{lemma}
\label{le:justone}
For any $i\in \Z_{>0}$ there exists $j\in \Z_{>0}$ such that $\phi(Z_i) \subset W_j$.
\end{lemma}
\begin{proof}
Observe that  $\bar \phi$ is a continuous map under  topologies of $Z$ and $V$.  Further any closed in relative topology subset $A_i$ of  $V_i$ is a closed subset of $V$.  Therefore $\bar \phi^{-1}(A_i)$ is a closed subset of $Z$.   Consequently, $\bar \phi^{-1}(A_i)\cap Z_j$ is also a closed subset of $Z$, and thus  is a closed subset of $Z_j$ in the relative topology inherited from $Z$.  Note also that $Z_j=\bigcup_{i}(\bar \phi^{-1}(W_i) \cap Z_j)$ which is a union of disjoint closed sets.  Since $Z_j$ is an irreducible closed set, we conclude that all but one term in the union is empty.
\end{proof}
We next show that $\phi$ is continuous in the topologies of $Y$ and $Z$
\begin{lemma}
\label{le:cont}
$\phi$ is  a continuous map.
\end{lemma}
\begin{proof}%

       Now let $S$ be an open subset of $V$ such that  $S_W=S \cap W$ is an open subset of $W$. In this case $S \cap V_j=S_j$ is a  non-empty subset of $V_j$.   Also $\phi^{-1}(S)$ is an open subset of $Z$.  It remains to show that $\phi^{-1}(S) \cap Z_i \not =\emptyset $ for any $i$.    Fix an $i$ and let $\phi(Z_i) \subset V_j$,  then by Definition \ref{def:subvarieties} we have that $\phi^{-1}(S_j)\cap Z_j \not = \emptyset$. Hence $\phi$ is continuous in the topologies of $Y$ and $W$.
\end{proof}
The final part we need to show that $\phi$ is a morphism of ringed space is to show that $\phi^{\#}$ is a homomorphism of rings.
\begin{lemma}
\label{le:homrings}
For any $O_W$ an open subset of $W$ we have that $\phi^{\#}:\calF_W(O_W) \rightarrow \calF_{Y}(\phi^{-1}(O_W))$ is a morphism of rings.
\end{lemma}
\begin{proof}
Since $\bar \phi$ is a composition morphism of ringed spaces, if $f \in \calF(O)$ for some open subset $O$ of $V$,   we have that $f\circ \phi \in \calF(\bar \phi^{-1}(O))$ and the map $f \mapsto f \circ \bar \phi$ is a homomorphism of the function rings.  Now if $O_W$ is an open subset of $W$ and  $f \in \calF(O_W)$, then $f \circ \phi$ is defined on $\phi^{-1}(O_W)$, an open subset of $Y$ by the discussion above.  Further, if $f \in \calI(W)$, then $f \circ \phi$ is identically zero on $\phi^{-1}(O_W)$ and thus is in $\calI(Y)$ by Remark \ref{rem:I}.
\end{proof}%
\begin{remark}
As in  Remark \ref{rem:justone}, we can apply the proposition above to the case when $Y=Z_i$ and conclude that $\bar \phi$ restricted to each $Z_i$ is a composition morphism of ringed spaces.
\end{remark}
From Lemmas \ref{le:justone} -- \ref{le:homrings} one can easily derive the following corollary.
\begin{corollary}
\label{cor:pbp}
In the notation above we have that
\begin{itemize}
\item if $T$ is a topologically irreducible ringed space of functions into $k$, then there exists a morphism from $T$ to $W$ as unions of  ringed subspaces if and only if there exists a composition morphism from $T$ to $V_i$ for some $i$ as ringed spaces of functions into $k$;
 \item $Y \cong W$  if and only if for each $i$ there exists a $j$ such that $Z_i \cong V_j$ as ringed spaces of functions into $k$ using a composition isomorphism; and
 \item $Y$ can be embedded into $W$ if and only if each $Z_i$ can be embedded into a distinct  $V_j$ as ringed spaces of functions into $k$ using a composition morphism.
 \end{itemize}
\end{corollary}
We now add some logical conditions on our unions of  ringed subspaces
to satisfy the conditions of Theorems \ref{thm:int1} and
\ref{thm:int2}, and prove, as the result, the main result of the
present section.

\begin{theorem}%
\label{prop:satisfy1}
Let $Z$ be of $\sZ$-class, and for each positive integer $i$, let $Z_i
\subseteq Z$ be closed and irreducible in the topology of $Z$.  Suppose that the following
additional conditions are also satisfied:
\begin{enumerate}
\item $Z$ is computable as a set of points;
\item $R(Z), K(Z)$ are computable as rings;
\item $k$ is computable as a subset of $R(Z)$ and $K(Z)$;
\item there exists a computable function $\chi : Z \longrightarrow \Z_{>0}$ such that given an element of $Z$ it outputs $i$ if and only if the element is in $Z_i$;
\item For each $X \subset \Z_{>0}$ and for each open subset $O$ of $Z$, there exist   functions $\phi_{O,X}$ and $\psi_{O,X}$ computable from the characteristic function of $X$ identifying  elements of $K(Z)$ that are elements of $\calF(O\cap Y)$ and $\calI(O \cap Y)$ respectively for $Y = \bigcup_{i \in X}Z_i$;
\item for any indices  $i \not = j$, there is no injective composition morphism from $Z_j$ to $Z_i$;
\end{enumerate}

Let $\calC$ be the class of ringed spaces of functions formed by
unions of  $Z_i$ under the topology described above.  Let
$\calA_i=Z_i$ and for a set $X \subset \omega$ let $\calA_X = \cup_{i
  \in X}Z_i$. Let $\calA = \emptyset$.  Then $\mathcal{C}$ admits
arbitrary degrees of isomorphism types, and admits isomorphism types
without degree.
\end{theorem}%

\begin{proof}%

Given our assumptions, in particular the assumption that subsets are disjoint, it is clear that $\mathcal{A}_X$ is Turing and enumeration reducible to
$X$.  Now let $\calB= \bigcup_{i\in I}Z_i, I \subset \Z_{>0} $ be such that $\calB \cong \calA_X$.
In this case by Corollary \ref{cor:pbp} the constituent sets of $\calB$ and $\calA$ are pairwise
isomorphic and given our assumptions are the same, and therefore $A_i$ is embeddable into $\calB$ only if $A_i$ is one of the
constituent subsets of $\calB$.  Whether or not this is the case is a computable  procedure in $\calB$
(requiring examination of points in $\calB$) by our assumptions.  Then $\mathcal{C}$ satisfies the conditions of Theorems \ref{thm:int1} and \ref{thm:int2}.
\end{proof}%

We now construct some examples of ringed spaces described above.  From the assumptions we made on our ringed spaces it would seem that
varieties are the most natural building blocks for our examples.
\begin{proposition}
\label{prop:computablevariety}
If $V$ is variety over a countable algebraically closed field $k$ presented by a finite set of equations, then as a ringed space, $V$ is of $\calZ$-class and is computable as a ringed space in our ringed space language uniformly form its equations.
\end{proposition}
\begin{proof}
We check the properties of the ringed spaces of $\calZ$-class:
\begin{enumerate}
\item The topological space $V$ is countable:  since $k$ is countable, the variety will have countably many points,
\item The number of open/closed subsets of $Z$ is countable.
\item $\calF(Z)=R(Z)$ is a countable ring of functions from $Z$ to some field $k$.  The ring $R(Z)$ is an integral domain and is also a $k$-algebra.  The fraction field of $R(Z)$ is denoted by $K(Z)$.
 \item For any open $U \subseteq Z$ we have that  $\calF(U) \subset
   K(Z)$ and $\calF(U)$ contains functions of $K(Z)$ defined on $U$ (but not necessarily all the functions defined on $U$).  A function $f \in K(Z)$ is defined on $U$ if $f = \frac{f_1}{f_2}$, with $f_1, f_2 \in R(Z)$ and $f_2(z) \not = 0$ for any $z \in U$.
\item  For  $U \subseteq V$  open sets, we have that   $\calF(V)\subseteq \calF(U)$ and $\rho_{U,V}$ is defined to be the inclusion map.
\item For any collection $A$ of open subsets of $Z$ we have that $\bigcap_{U \in A}\calF(U)=\calF(\bigcup_{U \in A}U)$.
 \end{enumerate}
  Since $k$ is countable, the variety will have countably many points, and each closed set will correspond to a finite system of polynomial equations over $k$ making the cardinality of the class of closed set countable.  If the rings corresponding to an open set consist of rational functions defined on the set, then every ring is an integral domain contained in the function field of the variety and is also countable.  Further, if $k$ is computable, then the set of points of the variety is uniformly computable from its equations, and given polynomial equations defining a closed set, the membership in the set is also computable. By Proposition \ref{prop:idealcomput} and Lemmas \ref{le:closedset}, any ring corresponding to an open subset of the variety is computable from the equations defining the complement of the given open set.  Finally, since the open sets are dense in Zariski topology,  any irreducible subvariety of a variety satisfies the zero-set condition.
\end{proof}
We should also note that any rational morphism of varieties will be a composition morphism, but not every rational morphism is, \emph{a priori}, continuous under the topology we set up for the union of subvarieties.  We remind the reader that we weaken the topology to contain only those non-empty open sets which intersect every subvariety in the union.  However, if we are in the situation where there are no rational morphisms between between subvarieties in two unions, we can conclude that there are no morphisms between the unions.
In our first example, we will consider a set open only if it is co-finite.   Note that this is a weakening of Zariski topology.  Below
we show that  weakening of the topology does not alter the status of the ringed
space.

\begin{lemma}[Weaker Topology]
\label{le:weaker}
Let $(Z,O_Z)$ be a ringed space under some topology $T$.  If $T_w$ is  a topology on $Z$ weaker than $T$, then $(Z,\tilde O_Z)$, where $\tilde O_Z$ contains only the rings corresponding to the open sets in $T_w$, is also a ringed space.   Furthermore, if $(Z, O_Z)$ satisfied the requirements for $\calZ$ class, then $(Z, \tilde O_Z)$ will also satisfy these requirements.
\end{lemma}
\begin{proof}
It is enough to note that every open set of $T_w$ is an open set of $T$ and every open covering of $T_w$ is also an open covering of $T$, while the rings corresponding to open sets common to both topologies are the same.
\end{proof}
\begin{remark}
\label{rem:co-finite}
If the new topology is co-finite topology, then any ringed space composition automorphism $\phi$ which is finitely-many-to-one remains a ringed space automorphism since for any such map $\phi$ the inverse image of a co-finite set is also co-finite.
\end{remark}

We are now ready to describe a general example where the underlying space will be an algebraic variety, though we will weaken its topology.
\begin{example}
Let $V$ be a  variety over an algebraically  closed field $k$,  let $Z = V$  as a set with a topology where all
the non-empty open sets are co-finite.  Let $R(Z)$ be the  coordinate
ring of the variety, computable by Proposition \ref{prop:idealcomput}
and set $\calF(U)$ for any open $U$ to be defined as usual as the ring of rational functions from $k(V)$ defined on $U$.  Let $\{Z_i\}_{i \in \omega}$ be a sequence of disjoint infinite subvarieties of $V$ such that there is a recursive procedure computing a set of polynomial equations defining $Z_i$ from $i$. We will use the rational morphisms of varieties satisfying Definition \ref{def:subvarieties} as morphisms between our union ringed subspaces.

    If $f \in k(V)$, then by Lemma \ref{le:infnum},   unless $f \equiv 0$ in $R$, we have that $f =0$ on finitely many $Z_i$ only.   Hence, if $Y = \bigcup_{i \in X}Z_i$ for some infinite set $X$ of positive integers, then a non-zero function in $R$ is not zero on $Y$.  If $X$ is a finite set, then  one can effectively determine whether $f$ is 0 in $Y$ by Lemma \ref{le:closedset1}.  Hence for any $Y= \bigcup_{i \in X} Z_i$, we have that $\calI(Y)$ is computable and therefore, $\calF_Y(U_Y)$ is computable from $X$.  The final assumption we need concerns possible morphisms between the subvarieties in the sequence.  We  we assume that that there is no non-constant morphism between any pair $(Z_i, Z_j)$ for $i \not = j$.   Given the last assumption we can conclude that $(Z, O_Z), \{Z_i\}$ now satisfy Proposition \ref{prop:satisfy1}.
\end{example}
We need the following lemma concerning non-singular curves to present more concrete examples of ringed spaces described above.

\begin{lemma}
\label{le:genus}
Let $C_1, C_2$ be two non-singular irreducible affine or projective curves over an algebraically closed field $k$.  Suppose, there exists an injective rational morphism
\[
\phi: C_1 \longrightarrow C_2.
\]
In this case the function fields of $C_1$ and $C_2$ are isomorphic and thus the curves have the same genus.
\end{lemma}
\begin{proof}
Since morphism is injective, it is non-constant and therefore it is onto.  (See [\cite{Sil1}, II, Theorem 2.5].)
\end{proof}
We now ready for our more concrete examples.  We start with a specific instance of the pattern described above.
\begin{example}%
Let  $Z$ be the elliptic surface  $y^2=x^3+x+z$ over $\tilde \Q$ -- the algebraic closure of $\Q$.   Note that all but two values of $z$ produce a non-singular curve, i.e. an elliptic curve.  (For a
Weierstrass equation of the form $y^2=x^3+Ax+B$, the discriminant is $\triangle=-16(4A^3+27B^2)$. See Remark 1.3, pages 49--50 of
\cite{Sil1}.) We can also make sure that no two such curves are isomorphic over $\tilde \Q$  by arranging for the different $j$-invariants ($j(E_{i})=\frac{-1728(4^{3})}{\triangle (E)}$).   Let $n_i$ be a computable sequence of positive natural numbers such that the sequence of curves $E_i = \{y^2=x^3 +x +n_i\}$ contains pairwise non-isomorphic (under rational morphisms) elliptic curves only.  Observe that the curves will be disjoint.   Let $Z_i=E_i$ and note that all the conditions of Proposition \ref{prop:satisfy1} are now satisfied as above.
\end{example}%
The next example has a somewhat different flavor and we need the following lemma to construct it.
\begin{lemma}
\label{le:genus1}
The  genus of the curve $y^{d} = (x+1)(x+2)$ over a field of constants $k$ of characteristic 0, where $d$ is an odd positive integer, is equal to $(d+1)/2$.
\end{lemma}
\begin{proof}
We consider the field extension $k(y,x)/k(x)$ and determine all the primes ramified in this extension.  Since $d$ is odd, the pole of $x$ is completely ramified in this extension. Examining the determinant of the extension yields that the only other ramified primes are the zeros of $x+1$ and $x+2$ and they are ramified completely.  Therefore the degree of the different is $3(d-1)$.  Using a version of Riemann-Hurwitz formula (see \cite{Fried}, Remark 3.6.2(c)), we now conclude that the genus of $k(x,y)$ is equal to $1 -d + 3(d-1)/2=\frac{d+1}{2}.$
\end{proof}

\begin{example}
\label{ex:genus}
Let $k$ be any computable algebraically closed field of characteristic 0, let  $Z=\A^3_k$ under co-finite topology, i.e. let a subset of $Z$ be open if it is co-finite
Let  $Z_n = \{(x,y,n): y^{d_n}=(x+1)(x+2)\}$, where $d_n$ is an increasing sequence.  The curve $y^{d_n}=(x+1)(x+2)$ is a non-singular curve (by a direct calculation) of genus $(d_n+1)/2$ by Lemma \ref{le:genus1}.)  Observe that the curves are pairwise disjoint and by Lemma \ref{le:genus} no curve in the sequence has a rational embedding into another curve.  Thus, by the discussion above, the sequence $\{Z_i\}$ satisfies the assumptions of Theorem \ref{prop:satisfy1}.

\end{example}

\section{The Class of Unions of Ringed Spaces}\label{secunionvar}
\setcounter{equation}{0}
In this section we will consider  ringed spaces which are unions of {\it arbitrary} ringed spaces, not necessarily subsets of a single ringed space.

Let $W = \bigcup_{j \in \Z_{>0}}V_j$, where each $V_j$ is a $\calZ$-class ringed space of functions into a countable field $k$, satisfying the zero-set condition, and for $m \not = j$ we have that $V_m \cap V_j =\emptyset$.  We use open subsets of all $V_j$ to construct a  topology on $W$: an open subset $U$ of $W$ will be of the form $U=\bigcup_{j \in \omega} U_j,$ where for all  $j$ we have that  $U_j$ is a non-empty  open subset of $V_j$ and for all but finitely many $j$ we have that $U_j = V_j$.  Any union and any finite intersection of sets of this form is again a set of this form, and therefore we  defined a topology on $W$.

\begin{proposition}[Ringed Space of a Union of Ringed Spaces]%
\label{prop:unions} Let $\calF_j$ be the ring of  functions on $U_j$,
as defined in the previous paragraph.  Now we let $f \in \calF(U)$ be a function from $W$ to $k$ such that for all  $U_j$ we have that $f|_{U_j} \in \calF(U_j)$ and for some $c \in k$ for all but finitely many $j$ we also have that $f|_{U_j} \equiv c$.    We claim that $\{\calF(U): U \mbox{ an open subset of } W\}$ is a sheaf of rings and $W$ is a ringed space.
 \end{proposition}%

 \begin{proof}%
 We check that the following conditions hold.
  \begin{enumerate}%
  \item {\it For any open $U \subset W$, we have that $\calF(U)$ is a ring.}   This clearly holds by construction.
\item For every inclusion $S_1 \subseteq S_2$ of open subsets of $W$  there exists a morphism of rings $\rho_{1,2} : {\calF}(S_2) \rightarrow {\calF}(S_1)$: in our case $\rho_{1,2}$ can be taken to be the inclusion map.
 \item {\it $\rho_{U,U}$ is the identity map ${\calF}(U) \longrightarrow   {\calF}(U)$}: this is true if we let $\rho$ be the inclusion map.
\item {\it If $W \subseteq V \subseteq U$ are three open sets, then $\rho_{U,W} = \rho_{V,W}\circ\rho_{U,V}$}:  this is true again due to the fact that $\rho_{U,W}, \rho_{V,W}, \rho_{U,V}$ are all inclusion maps.
\item {\it  If $U$ is an open set, $\{B_i, i \in A \subset \Z_{>0}\}$  is an open covering of $U$, and  $s \in {\calF}(U)$ is an element such that $\rho_{U,B_i}(s)=0$ for all $i$, then $s=0$.}  This follows from the fact that $\rho_{U,B_i}$ is an embedding.
\item  {\it If $U$ is an open set, $\{B_i, i \in A \subseteq \Z_{>0}\}$ is an open covering of  $U,$ and for each pair $i, j$, there exist $s_i \in {\calF}(B_i)$ and $s_j \in {\calF}(B_j)$ such  that
    \[
    \rho_{B_i,B_i\cap B_j}(s_i)= \rho_{B_j,B_i \cap B_j}(s_j),
    \]
    then there exists $s \in {\calF}(U)$ such that for each $i$, we have that $\rho_{U,B_i}(s)= s_i$. } Let $G(Z)$ be defined as above as a ring generated by all the rings $\calF(U), U \subseteq W$ open.  If $\rho_{B_i,B_i\cap B_j}(s_i)= \rho_{B_j,B_i \cap B_j}(s_j)$, then $s_i = s_j \in G(Z)$, and therefore
\[
    s_j \in \calF(B_j) \cap \calF(B_i).
\]
Fixing $i$ and varying $j$ we obtain that $s_i \in \bigcap_j \calF(B_j)$ and thus
\[
s_i|_{V_m}  \in  \bigcap_j \calF(B_j \cap V_m).
\]
Using a property of the $\calZ$-class ringed spaces we also conclude that
\[
s_i|_{V_m}  \in   \calF(\bigcup B_j \cap V_m)=\calF(U_m).
\]
Thus, $s_i \in \calF(U)$.
    \end{enumerate}
\end{proof}
\begin{definition}[Morphisms of Unions of Ringed Spaces]
\label{def:morunion}
Let $$W_1=\bigcup_{i \in \Z_{>0}}V_{i,1}, W_2=\bigcup_{i \in \Z_{>0}}V_{i,2}$$ be two ringed spaces which are unions of ringed spaces under the topology described above.  In this case  a map $\phi : W_1 \longrightarrow W_2$ is a morphism of unions of ringed spaces if for every $i$ there exists  $j$ such that  $\phi_i= \phi|_{V_{i,1}} :V_{i,1} \longrightarrow V_{j,2}$ is a composition morphism of  ringed spaces and for any $j$ we have that  $\phi^{-1}(V_{j,2}) \cap V_{i,1}= \emptyset$ for all but finitely many $i$. If $\phi$ is injection, we call the corresponding morphism an embedding and if  $\phi$ is a bijection and its inverse is also a morphism of  unions of ringed spaces, then $\phi$ will be called an isomorphism of unions of ringed spaces.
\end{definition}
\begin{lemma}
\label{le:welldefined}
Morphism of unions of varieties as defined above is a morphism of ringed spaces.
\end{lemma}
\begin{proof}
Let $W_1, V_{i,1}, V_{j,2}, W_2, \phi$ be as in Definition \ref{def:morunion}.  If $U_2 \subset W_2$ is a non-empty open set, then $U_2=\bigcup_{j \in \Z_{>0}}U_{j,2}$, where $U_{j,2}$ is a non-empty open subset of $V_{j,2}$ in (relative topology).  Observe that
\[
\phi^{-1}(U_2) = \bigcup_{j \in \Z_{>0}}\phi^{-1}(U_{j,2})= \bigcup_{j \in \Z_{>0}}\bigcup_{i \in I_j}\phi_i^{-1}(U_{j,2}),
\]
where for each $U_{j,2}$ we let $I_j=\{i \in \Z_{>0}: \Im (\phi_i) \subseteq V_{j,2}\}$.  Observe that $\bigcup_{j \in \Z_{>0}} I_j =\Z_{j>0}$ and $\phi_i^{-1}(U_{j,2})$ is an open (in relative topology) set in $V_{i,1}$, and therefore $\phi$ is a continuous map.

Further, for any open set $U_2 \subset V_2$ and $f \in \calF(U_2)$ we have a ring homomorphism
\[
\phi^{\#} : \calF(U_2) \longrightarrow \calF(\phi^{-1}(U_2)),
\]
defined by $f \mapsto f\circ \phi$.   Indeed,  let $U_1=\phi^{-1}(U_2)$, let $U_1 \cap V_{i,1}=U_{i,1}$, let $\phi_i(U_{i,1}) \subseteq V_j$, let $U_{j,2} =U_2 \cap V_{j,2}$  and note that $U_{i,1} = \phi_i^{-1}(U_{j,2})$.    Observe now that $f\circ \phi|_{U_{i,1}}= f\circ \phi_i \in \calF(\phi^{-1}(U_{j,2}))$  as required and if $f|_{U_{j,2}}=c \in k$ for all but finitely many $j$, then $f\circ \phi|_{U_{i,1}}=c$ for all but finitely many $i$, since only finitely many sets $U_{i,1}$ are mapped into any given $U_{j,2}$.
\end{proof}

Given our definition we immediately obtain the following propositions.

\begin{proposition}
  Let $V, W, V_{i,1}, V_{j,2}, k$ be as above.   Then $W_1$ is embeddable into $W_2$ if and only if for each $i$ there exists a $j$ such that $V_{i,1}$ is embeddable into $V_{j,2}$.  Further, $W_1 \cong W_2$ if and only for each $i$ there exists a $j$ such that $V_{i,1} \cong V_{j,2}$ as ringed spaces via composition morphism.
\end{proposition}

\begin{proposition}%
\label{prop:satisfy11}
Let $W=\bigcup_{i=1}^{\infty} V_i$, be a union of ringed spaces  as in Proposition \ref{prop:unions} with additional conditions
conditions:
\begin{enumerate}
\item the field $k$ is computable;
\item the set of all points of $W$ is computable
\item there exists a function $f(i,x)$ which given an index $i$ and an element $x \in W$ determines whether $x \in V_i$;
\item given $V_j, V_i$, $i \not = j$, there is no injective ringed space morphism  from $V_j$ to $V_i$;
\item  $V_i$ are uniformly computable as  ringed spaces, i.e. there is a computable function $\psi(i,j)$ determining  the whether the $j$-th statement in the atomic diagram of $V_i$ is true.
\end{enumerate}

Let $\calC$ be the class of ringed spaces formed by unions of ringed spaces contained in $W$.  Let $\calA_i=V_i$
and for a set $X \subset \omega$ let $\calA_X = \cup_{i \in
  X}V_i$. Let $\calA = \emptyset$.  Under these assumptions
$\mathcal{C}$ admits arbitrary degrees of isomorphism types, and
admits isomorphism types without degree.
\end{proposition}%

\begin{proof}%

Given our assumptions it is clear that $\mathcal{A}_X$ is Turing and enumeration reducible to
$X$.  Now let $\calB$ be a union of  some ringed spaces from $W$ such that $\calB \cong \calA_X$.
Then by Corollary \ref{cor:pbp} the constituent ringed spaces of $B$ and $\calA$ are pairwise
isomorphic, and therefore $A_i$ is embeddable into $\calB$ only if $A_i$ is  one of the
constituent varieties of $\calB$.  Whether or not this is the case is
an r.e.\ condition in $\calB$
(requiring examination of the coordinates of points in $\calB$) by our
assumptions.  Thus, the conditions of Theorems \ref{thm:int1} and
\ref{thm:int2} are satisfied.

\end{proof}%

We have established the following.

\begin{theorem}
\label{thm:unionvar}
 In the notation above let $\calC$ be the class of ringed spaces of functions formed by unions of ringed spaces taken from the collection $\{V_i\}$.   Then for any Turing degree $\mathbf{d}$, there is a member of $\calC$ whose isomorphism type has
Turing degree $\mathbf{d}$.  There is also a member of $\calC$ whose
isomorphism type has no Turing degree.
\end{theorem}
We now give specific examples of varieties illustrating Theorem \ref{thm:unionvar}.   First we note that from Proposition \ref{prop:idealcomput} and Lemma \ref{le:closedset} we can derive the following assertion.
\begin{proposition}
\label{prop:curveequation}
If $C$ is a plane irreducible curve over an algebraically closed
computable field $k$ given by a single equation, then $C$ under the
canonical ringed space structure is computable as a ringed space in
our ringed space language uniformly in its equation.
\end{proposition}

  First we reconsider the curve examples from above with one significant difference: the curves will keep their natural ringed space structure.

\begin{example}%
Consider an elliptic surface $y^2=x^3+zx+1$ over $\tilde \Q$ -- the algebraic closure of $\Q$.  All
but three values of $z$ produce a non-singular curve, i.e. an elliptic
curve.  By considering an increasing sequence of values for $z$, as
above, we can make sure that the corresponding curves have different
$j$-invariants and thus are not isomorphic over $\tilde \Q$.  Let
$(n_i)_{i \in \omega}$ be such a sequence which is also computable.  Let $V_i : y_i =x_i^3+n_ix +1$ with the indigenous ringed space structure,
and observe that  by Proposition \ref{prop:curveequation}, all the conditions of Proposition \ref{prop:satisfy11} are satisfied.
\end{example}%

\begin{example}
If  $\{V_i\}_{i \in \omega}$ is a sequence of plane irreducible curves of genus $g_i$, with $g_i \not = g_j$ for $i \not = j$, under the natural ringed space structure, such that the sequence of their equations is computable, then all the conditions of Proposition \ref{prop:satisfy11} are satisfied.  (For an example of a specific sequence we can us the sequence constructed in Example \ref{ex:genus}.)
\end{example}

We can specialize our examples  above to the case of Abelian varieties.  Here the constituent parts of the union will be an Abelian variety and the morphisms considered will restrict to an isogeny for any constituent part of the union.  The example below uses products of varieties along the lines of an idea in \cite{Shafarevich}.

\begin{example}\label{needapp1}
Let $\{a_{i,j}\}_{i,j \in \omega}$ be a double indexed computable collection of pairwise non-isogenous Abelian varieties over an algebraically closed field $k$,  and let $V_i = \prod\limits_{j \leq i}a_{i,j}$.   We claim that all the conditions of Proposition \ref{prop:satisfy11} are satisfied.  Indeed the only condition we need to check is the absence of non-zero isogenies between $V_i$ and $V_j$ for $i\not=j$.

Observe that  there is an  injective isogeny from $a_{r,m}$ into $V_i$  if and only if $r=i$ and $m \leq i$.  Indeed for $j \leq i$ let $\pi_{i,j}$ be the projection from $V_i$ onto $a_{i,j}$.  Now let $\phi : g_{r,m} \rightarrow V_i$ be an injective isogeny.  In this case $\pi_{i,j} \circ \phi : a_{r,m} \longrightarrow a_{i,j}$ and given our assumptions, unless $(i,j)=(r,m)$, this isogeny must be the zero isogeny.  Now since $\phi$ is not constant,  for some $j \leq i$ we have that the composition is not constant, and thus the assertion follows.

  Note further that there exists a non-constant isogeny from $V_i$ to $V_j$ if and only if $i=j$.   Otherwise for some $k \leq i$ we have that there exists a non-constant isogeny from $a_{i,k}$ to $V_j$. Indeed, let $(\bar x_1,\ldots \bar x_i)$ be the  coordinates of points of $V_i$, and $(\bar y_1,\ldots, \bar y_j)$ be the coordinates for $V_j$.   Let $(\bar \phi_1,\ldots,\bar \phi_i)$ be the isogeny from $V_i$ to $V_j$.    Consider $\bar \phi_1(\bar x_1,\ldots, \bar x_i)$ as we fix values of $\bar x_2, \ldots, \bar x_i$ and run through all the possible values of $\bar x_1$.  Either for some tuple of fixed values of $(\bar x_2,\ldots, \bar x_i)$ the set $\bar \phi_1(\bar x_1,\ldots, \bar x_i)$ contains more than one value (i.e. zero-value), or this is not the case.  If this is not the case, then for some $\bar a_1 \in \PP^2$ we have that
  \[
  \bar \phi_1(\bar x_1,\ldots,x_i) =\psi_1(\bar a_1,\bar x_2, \ldots, \bar x_i)=\bar \phi_{1,2}(\bar x_2\ldots, \bar x_i).
\]
Now we repeat our search for distinct values but this time for $\bar \phi_{1,2}$ running through all the values of $\bar x_2$ while keeping $\bar x_3,\ldots,\bar x_i$ fixed.  Either one of these sets has more than one value or
\[
\bar \phi_1(\bar x_1,\ldots,\bar x_i) =\psi_1(\bar b_1,\bar b_2, \bar x_3,\ldots, \bar x_i)=\bar \phi_{1,3}(\bar x_3\ldots, \bar x_i).
\]
Repeating this procedure at most $i$-times we  will either find that for some $r \leq i$ and some set of values $(\bar b_1,\ldots,b_{r-1},\bar b_{r+1},\ldots, \bar b_i)$ the set $\{\bar \phi_1(\bar b_1,\ldots,\bar b_{r-1},\bar x_r, \bar b_{r+1},\ldots, \bar b_i)\}$ has more than one value as $\bar x_r$ takes all possible value, or in fact  $\bar \phi_1$ is constant.  Since All $\bar \phi_m, m=1,\ldots,j$ cannot be constant by assumption, for some $m \leq j$ we have that $$\{\bar \phi_m(\bar b_1,\ldots,b_{r-1},\bar x_r, \bar b_{r+1},\ldots, \bar b_i)\}$$ is not constant.  Thus we have a non-constant isogeny from $g_{i,r}$ into $V_j$ which is impossible unless $i=j$.

Finally we note that a computable collection of pairwise non-isogenous curves is constructed in the Appendix.
\end{example}%

\section{The Ringed Spaces of ``Elements"} \label{secschemes}

In this section we discuss a different family of ringed space for which Theorems \ref{thm:int1} and
\ref{thm:int2} hold: the schemes.  The category of schemes admits such
diversity that the rich structure guaranteed by these theorems may arise
in a number of different ways.  Correspondingly, we give several
constructions in the present section, some very familiar, some less so.  We start with a definition of a class of ringed spaces which contains the $\calZ$-class.
\begin{definition}[Ringed space of ``elements'']
\label{def:subspace}
Let $Z$ be a countable irreducible topological space with a countable topology. Let $(Z, O_Z)$ be a sheaf of rings on $Z$ satisfying the following conditions:
 \begin{enumerate}
 \item There exists a countable ring $R(Z)$ such that for any open $U \subseteq Z$ we have that  $\calF(U) \subset R(Z)$.
 \item If $U \subseteq V$ are non-empty open sets, then $\calF(V)\subseteq \calF(U)$ and the map $\rho_{V,U}$ is the inclusion map.
 \item For any finite collection of open sets $U_i$  we have that $\calF(\bigcup_i U_i) = \bigcap_i \calF(U_i)$.
 \end{enumerate}
 \end{definition}
 Next we discuss a property of subsets of a ringed space.  The proposition below is quite similar to Proposition  \ref{prop:well-def}.
 \begin{proposition}
 If $Y \subset Z$ is any subset of $Z$, irreducible under the relative topology, then  $Y$ also has a sheaf structure on it with an inclusion map morphism into $(Z, O_Z)$.  Furthermore, this sheaf structure satisfies Definition \ref{def:subspace}.
\end{proposition}
\begin{proof}
Given an open set $U_Y$ of $Y$, let $\calF(U_Y) = \bigcup_{U_Y \subset U}\calF(U)$, where the union is taken over all open sets $U$ of $Z$ containing $U_Y$.  It is not hard to see that this union is actually a sub-ring of $R(Z)$.  Indeed let $U$ be an open subset of $Z$ such that $U \cap Y \supseteq U_Y$ (at least one such set exists by definition of relative topology).  Since $\calF(Y)$ is a ring, it contains $1_{R(Z)}, 0_{R(Z)}$ (multiplicative and additive identity from $R(Z)$ respectively), and therefore $1_{R(Z)}, 0_{R(Z)} \in \calF(U) \subset \calF(U_Y)$.  Further, if $a, b \in \calF(U_Y)$, then for some open $U_a, U_b \subset Z$ we have that $U_Y \subset U_a \cap U_b$ and $a \in \calF(U_a), b \in \calF(U_b)$.  Thus, $\calF(U_a \cap U_b)$ contains both $a$ and $b$ and therefore $ab$ and $a+b$.  Consequently, $ab$ and $a+b$ are both in $\calF(U_Y)$.
It is also not hard to see from the definition of $\calF(U_Y)$ and $\calF(V_Y)$ for open subsets $U_Y$ and $V_Y$ of $Y$ that if $U_Y \subseteq V_Y$ are open subsets of $Y$, then as above  we have $\calF(V_Y)\subseteq \calF(U_Y)$ and the map $\rho_{V_YU_Y}$ can be defined to  be the inclusion map, again as above.

Next we show that if $U_Y$ and $W_Y$ are open subsets of $Y$ then
\[
\calF(U_Y \cup W_Y) = \calF(U_Y) \cap \calF(W_Y).
\]
    Indeed,
\[
 \calF(U_Y) \cap \calF(W_Y)= \bigcup_{U_Y \subseteq U}\calF(U) \cap \bigcup_{W_Y \subseteq W}\calF(W)=\bigcup_{U_Y \subseteq U, W_Y \subseteq W}\calF(U) \cap \calF(W)=
 \]
 \[
\bigcup_{U_Y \subset U, W_Y \subseteq W}\calF(U \cup W)=\bigcup_{U_Y \cup W_Y \subseteq T}\calF(U \cup W)=\calF(U_Y \cup W_Y),
\]
where $T$ ranges over all open subsets of $Z$ containing $U_Y$ and $W_Y$, and the penultimate equality holds because for any open subsets $U, W$ of $Z$ with $U_Y \subseteq U$ and $W_Y \subseteq W$ we certainly have that $T=U\cup W$ is open and contains $U_Y \cup W_Y$, while for any open $T \subset Z$ with $U_Y \cup W_Y \subseteq T$, we  certainly have $U_Y  \subseteq T$ and $U_Y  \subseteq T$.

Now let $B_{Y,i}$ be an open covering of $U_Y$ and suppose that there exists an $s \in \calF(U_Y)$ such that for all $i$ we have that $\rho_{U_Y,B_{Y_i}}(s) =0_{\calF(B_{Y,i})}$.  In this case, since
\[
s=\rho_{U_Y,B_{Y_i}}(s) =0_{\calF(B_{Y,i})}=0_{R(Z)},
\]
we reach the desired conclusion that $s=0_{R(Z)}$.  Finally suppose that for all $i $ we have $s_i \in \calF(B_{Y,i})$ and for all $i,j$ it is the case $\rho_{B_{Y_{i}},B_{Y_{i}} \cap B_{Y_{j}}}(s_i) =\rho_{B_{Y_{j}},B_{Y_{i}} \cap B_{Y_{j}}}(s_j)$.  Observe that since $Y$ is irreducible in the relative topology, by Lemma \ref{le:intersec}, we have that $B_{Y_{i}} \cap B_{Y_{j}}$ is non-empty, and therefore we conclude that $s_j=s_i \in \calF(B_{Y_{i}}) \cap \calF(B_{Y_{j}})=\calF(B_{Y_{i}} \cup B_{Y_{j}})$.  Thus, keeping $i$ constant and varying $j$ we conclude that $s_i \in \bigcap_j\calF(B_{Y_j})$.  Applying the same argument as in Proposition \ref{prop:well-def}, we conclude that $s_i \in \calF(U_Y)$.

 Note further that the inclusion  map $i$ of $Y$ into $X$ produces a morphism of ringed spaces.  First of all, the inverse image of an open set  will be open, making this map is continuous.  Indeed, suppose $i^{-1}(O)=O_Y$, where $O$ is an open subset of $Z$  and $O_Y$ is an open subset of $Y$ and note that this means that $O_Y= O \cap Y$, making $O_Y$ open in the relative topology.  Note further, that $\calF(O) \subset \calF(O_Y)$ by construction, and thus we can define $i^{\#}$ to be the inclusion map of rings.
\end{proof}
To be able to apply Theorems \ref{thm:int1} and \ref{thm:int2}, we need to make one more assumption:
\begin{assumption}
\label{assum:irr}
Assume that  either {\it any} $Y \subset Z$ containing a certain subset is closed and irreducible under the relative topology,  or any infinite subset of $Z$ is irreducible under the relative topology.
\end{assumption}
\begin{remark}
A topology satisfying Assumption \ref{assum:irr} exists on any set.    To satisfy the first assumption one can let a set be open if it contains a certain subset of $Z$, possibly finite.  To satisfy the second assumption one can let the open sets be the co-finite sets.
\end{remark}

We now describe a class of schemes which admit arbitrary Turing
degrees, establishing the final point of the main result.

\begin{theorem}
\label{thm:subspace}
Assume $(Z, O_Z)$ is infinite, computable and satisfies Assumption
\ref{assum:irr}.   If every open set contains a certain subset, denote
this subset by $A$,  assume that it is computable and that its complement infinite. If there is no such subset,  let $A=\emptyset$. If $A$ is finite, let $B \subset Z$ be a computable subset of $Z\setminus A$ with an infinite complement.  If $A$ is infinite, let $B = \emptyset$.  Let $e: \Z_{>0} \longrightarrow Z\setminus (A\cup B)$ be a computable enumeration of  $Z\setminus \{A \cup B\}$  and for any $X \subset \Z_{>0}$ set
\[
Z_X=\{e(n): n \in X\} \cup B \cup A.
\]
Assume that $(Z_X, O_{Z_X}) \leq_T X$ and $(Z_X, O_{Z_X}) \leq_e X$, or in other words, for each open set $U$, the functions $\phi_{U,X}(f)$ determining membership in $\calF(U_X)=\calF(U \cap Z_X)$ and $\psi_{U,X}(f)$ listing the elements of $\calF(U_X)$ are computable from $X$. Consider the class  $\{Z_X: X\subseteq \Z_{>0}\}$ together with ringed space morphisms making the following diagram commute, where $i_M, i_N$ are the inclusion maps.
\begin{equation}
\label{eq:diag1}
{\xymatrix{
Z_N\ar[rr]^{i_N}&& Z\\
&Z_{M}\ar[ul]\ar[ur]^{i_M}&
}}
\end{equation}
Under the assumptions above the following statements are true.
\begin{enumerate}
\item There exists a ringed space morphism $i_{M,N}: Z_M
  \longrightarrow Z_N$ if and only if $M \subset N$.
\item Given $\mathcal B = Z_X$ for some $X$ and $Z_i \not =Z_j$, only one of which is embeddable into $\calB$, we can determine which of the two ringed spaces is embeddable into $\mathcal B$ by checking whether  $e(i)$ or $e(j)$ is contained in $\mathcal B$.
\item  The class of isomorphic copies of structures of the form $Z_S$ admits arbitrary degrees of isomorphism types,
  and admits isomorphism types without degree.

\end{enumerate}
\end{theorem}

\subsection{Examples Using Affine and Projective Schemes}

Our initial examples will deal with affine and projective schemes arising directly from ring spectra.
\begin{proposition}
\label{prop:schemes}
Let $K$ be a countable computable field, let $\{R_i, i \in \Z_{>0}\}$ be a collection of computable in $K$  local rings contained in $K$ such that
\begin{enumerate}
\item the fraction field of each $R_i$  is $K$,
\item $\pp_i$ is the maximal ideal of $R_i$, and each $\pp_i$ is computable in $R_i$.
\item for any $I \subsetneq \Z_{>0}$ and $R_I=\bigcap_{i \in I}R_i$ we have that $\Spec R_I =\{\pp_i \cap R_I, i \in I\} \cup \{(0)\}$,
\item for any $I \subseteq \Z$ we have that $R_I \leq _T I$ and $R_I \leq _e I$.
\end{enumerate}
Now let $Z=\{\pp_i, i \in \Z_{>0}\} \cup \{(0)\}$ under Zariski topology, i.e. a set $A= \{\pp_i, i \in I\}$ is closed if $A$ is the set of prime ideals containing an ideal $\mathfrak I$ of $R_J$ for some set of positive integers $J$ containing $I$. For any $X \subset \Z_{>0}$ set
\[
Z_X=\{\pp_{(2n+1)}: n \in X\} \cup \{\pp_{2n}: n \in \Z_{n>0}\} \cup \{(0)\}
\]
  Under these assumptions the above  described collection $\{(Z_X, O_{Z_X})\}$ satisfies Theorem \ref{thm:subspace}.
\end{proposition}
\begin{proof}
Below we check the  assumptions of Theorem \ref{thm:subspace}:
 \begin{enumerate}
 \item There exists a countable ring $R(Z)$ such that for any open $U \subseteq Z$ we have that  $\calF(U) \subset R(Z)$:   in our case the field $K$ contains all the rings.
 \item If $U \subseteq V$ are non-empty open sets, then $\calF(V)\subseteq \calF(U)$ and the map $\rho_{V,U}$ is the inclusion map: by definition we have that
  \[
  \calF(V) = \bigcap_{\pp \in V} R_{\pp} \subseteq \bigcap_{\pp \in U}R_{\pp}= \calF(U)
  \]
 \item For any collection of open sets $U_i$  we have that $\calF(\bigcup_i U_i) = \bigcap_i \calF(U_i)$: this part is true by construction.
 \item All the open sets contain the zero ideal.  If the rings above are Noetherian, then all the open sets are also co-finite.
 \item Given our assumptions, the ringed space corresponding to any $I \subset \Z$ is computable from $I$.
 \end{enumerate}

\end{proof}

 We now consider an example whose affine part is really a reconsideration of an example from \cite{CHS} where it was described in purely algebraic terms.  The projective part of the example is new.  However, the proofs are unchanged by the introduction of the projective part.
\begin{proposition}
\label{prop:projschemes}
Let $K$ be a computable product formula field, i.e. a number field or a finite extension of a rational function field over a computable field of constants.  Let $R_i, i \in \Z_{>0}$ be a collection of valuations rings of $K$ with valuations trivial on the constant field in the case $K$ is a function field.  Let $Z$ be the collection of all the valuation ideals under Zariski topology.   In this case, $(Z, O_Z)$ satisfies the conditions of Proposition \ref{prop:schemes}.
\end{proposition}
\begin{proof}
First of all we observe that if $K$ is a number field, then $(Z, O_Z)$, as well as  $(Z_X, O_{Z_X})$ for any $X\subseteq \Z_{>0}$ is an affine scheme and under Zariski topology  all the non-empty open sets contain the zero ideal and are co-finite. Therefore any infinite subset of $Z$ containing the zero ideal is irreducible under the relative topology.  Secondly, given a collection of prime ideals, from the characteristic function of the collection we can compute the characteristic function of the corresponding ring, i.e. the ring which is the intersection of the localization of  the ring integers of $K$ at the ideals in the collection. (See \cite{CHS} for more details.)  Finally we note that for a number field $R_{\Z_{> 0}}$ is the ring of integers of the field.

The situation is slightly different in the case $K$ is a function field.  If $K$ is a function field, then $(Z, O_Z)$ is a projective scheme, but for any $X \subsetneq \Z_{>0}$ we still have that $(Z_X, O_{Z_X})$ an affine scheme.  Further, $R_{\Z_{> 0}}$ is the constant field.  However, as in the case of number field, all the conditions of Proposition \ref{prop:schemes} are satisfied and we again refer the reader to \cite{CHS} for details.
\end{proof}

\begin{remark}
Using proposition Proposition \ref{prop:schemes} one can produce examples where $K$ is an infinite algebraic extension of $\Q$ or  an infinite algebraic non-constant extension of a rational function field.   However, to make sure that for all $I \subseteq Z$ we have that $R_I \leq_T I$  the set of prime ideals $Z$ should be selected so that only finitely many elements of $Z$ lie above a single rational prime.  We will not present details of such an example here as it is not sufficiently different from the example involving product formula fields.  Below we look at a case where  the transcendence degree over the constant field is greater than one.
\end{remark}

  We continue with a new set of  notation and assumptions.
\begin{notationassumptions}
\label{not:irr}\rule{5pt}{0pt}
\begin{itemize}
\item Let $k$ be an algebraically closed computable field.
\item Let $V$ be  variety over $k$ with $V \subseteq \A_k^{n+1}$.  We will assume that either  $V=\A_k^n$ or $V$ is the set of $k$-zeros of a polynomial of the form $F(x_1,\ldots,x_n,z)=\sum_{j=0}^{k}f_j(\bar x)z^j$, where $\bar x=(x_1,\ldots,x_n), f_j(\bar x) \in k[\bar x]$, and $f_k = 1$.   Further we will assume that for any $c \in k$ and any $j=1,\ldots,n$ we have that  $F(z,x_1,x_2,x_{j-1},c,x_j,\ldots,x_n)$ is irreducible in $k[z,x_1,x_2,x_{j-1},x_j,\ldots,x_n]$.
\item Let $k_j=k(x_1,\ldots,x_{j-1},x_{j+1},\ldots,x_n)$.
\item Let $\phi : \Z_{>0} \rightarrow k$ be a computable listing of $k$.
\item Let $R$ be the coordinate ring of $V$ and let $K$ be the fraction field of $R$.
\item Let $\bar R$ be the integral closure of $k[\bar x]$ in $K$.
\item Let $Z=\Spec R$ be the corresponding affine scheme.
\item  Let $I_1,\ldots, I_n$ be subsets of $\Z_{>0}$.
\item Let $P_{I_1,\ldots,I_n}$  be the following set
$$\{h \in k[\bar x]: \forall j =1,\ldots,n, \forall a \in \phi(I_j) \mbox{ we have } h \not \equiv 0 \mod (x_j-a) \mbox{ in } k[\bar x]\}.$$
\item Let $R_{I_1,\ldots,I_n} = \{g \in k(\bar x): g =\displaystyle \frac{g_1}{g_2}, g_1, g_2 \in k[\bar x], g_2 \in P_{I_1,\ldots,I_n}\}.$
\item Let $\bar R_{I_1,\ldots,I_n}$ be the integral closure of $R_{I_1,\ldots,I_n}$ in $K$.
\item Let $Z_{I_1,\ldots,I_n} = \Spec \bar R_{I_1,\ldots,I_n}$.
\item Let $K^G$ be the Galois closure of $K$ over $k(\bar x)$.
\item Let $\bar R^G$ be the integral closure of $\bar R$ in $K^G$.
\item Let $z=z_1,\ldots,z_k \in K^G$ be all the conjugates of $z$ over $k(\bar x)$.
\item Let $\bar R_j$ be the integral closure of $k_j[x_j]$ in $K$.
 \item Let $R_{a,j}$ be the localization of $k[\bar x]$ at $x_j-a$, and let $\bar R_{a,j}$ be the integral closure of $R_{a,j}$ in $K$.
\end{itemize}
\end{notationassumptions}

 Before we can proceed we need to review some basic facts from algebra.

\begin{lemma}
\label{le:basic}
Let $B$ be an integral domain with fraction field $F$.  Let $H$ be a finite extension of $F$, $M$ a finite extension of $H$,  and let $A_H$ the integral closure of $B$ in $H$ and let $A_M$ be the integral closure of $A_H$ in $M$.  We claim the following statements are true.
\begin{enumerate}
\item $A_M$ is the integral closure of $B$ in $M$.
\item If $a, b \in H$ and are integral over $B$, then $a+b$ and $ab$ are also integral over $B$.
\item If $\pp$ is a prime ideal of $B$, then there exists a prime ideal $\Pp$ of $A_H$ such that $\Pp \cap B=\pp$.
\end{enumerate}
\end{lemma}

\begin{proof}
Part 1 is Proposition 1.3, Section 1, Chapter VII of \cite{L2}.  Part
2 is Proposition 1.4, Section 1, Chapter VII of \cite{L2}.  Part 3 is Proposition 4.15, Section 4.4, Chapter 4 of \cite{Eisenbud}.
\end{proof}

Below we discuss some properties of $\bar R_{I_1,\ldots,I_n}$ and  $Z_{I_1,\ldots,I_n}$.
\begin{remark}
$R_{I_1,\ldots,I_n}$ is well defined since $P_{I_1,\ldots,I_n}$ is closed under multiplication.
\end{remark}

\begin{remark}
Since $k$ is algebraically closed, we can factor polynomials effectively in $$k[x_1,\ldots,x_n]$$ and over $$k(x_1,\ldots,x_n).$$ (See discussion of fields with splitting algorithms in Sections 17.1, 17.2 of \cite{Fried}.)
\end{remark}
\begin{proposition}
\label{prop:Turing1}
$R_{I_1,\ldots,I_n}$ and $\bar R_{I_1,\ldots,I_n}$ are both Turing
equivalent to the join of the $I_i$'s.
\end{proposition}
\begin{proof}
First of all we note that since $k$ is computable,  $k(\bar x)$ is computable and $k[\bar x]$ is computable as a subset of $k(\bar x)$.   At the same time, $K$ is also computable as a finite extension of $k(\bar x)$ with $k(\bar x)$ and $k[\bar x]$ being computable subsets of $K$. Further as remarked above we have an effective procedure for factoring  elements of $k[\bar x]$ over $k$. Consequently, given an element of $k(\bar x)$ we can write it  as a ratio of relatively prime polynomials and then using membership in $I_1,\ldots,I_n$ decide whether the element has any ``forbidden'' factors in the denominator, thus deciding whether the given rational function  is an element of $R_{I_1,\ldots,I_n}$.  Conversely, given the membership function in $R_{I_1,\ldots,I_n}$ we can determine if a $k$-element $a \in \phi(I_j)$ by checking whether
\begin{equation}
\label{eq:member}
\displaystyle \frac{1}{x_j -a} \in R_{I_1,\ldots,I_n}
\end{equation}
This argument shows that $R_{I_1,\ldots,I_n}$ is Turing equivalent to $\max(I_1,\ldots,I_n)$.  Next we show that $R_{I_1,\ldots,I_n}$ is Turing equivalent to $\bar R_{I_1,\ldots,I_n}$.  Since $R_{I_1,\ldots,I_n} = \bar R_{I_1,\ldots,I_n}\cap k[\bar x]$, we have that $R_{I_1,\ldots,I_n}$ is Turing reducible to $\bar R_{I_1,\ldots,I_n}$.  At the same time, using a power basis of $z$, given an element $y \in K$ we can determine the monic irreducible polynomial of $y$ over $k(\bar x)$ and using membership function of $R_{I_1,\ldots,I_n}$ determine whether $y$ belongs to $\bar R_{I_1,\ldots,I_n}$, the integral closure of $R_{I_1,\ldots,I_n}$ in $K$.
\end{proof}
In a similar fashion we can also establish the following equivalence.
\begin{lemma}
$R_{I_1,\ldots,I_n}$ and $\bar R_{I_1,\ldots,I_n}$ are both enumeration equivalent equivalent to $$\max(I_1,\ldots,I_n).$$
\end{lemma}

Next we need a couple of technical observations to be used below.
\begin{lemma}
\label{le:bigr}
For a pair of rings $R_1 \subset R_2$, if  $\Pp \subset R_2$ is a prime ideal, then $\Pp \cap R_1$ is also a prime ideal.
\end{lemma}
\begin{proof}
Suppose the lemma does not hold and therefore for some $z \in \Pp \cap R_1$ we have that $z =z_1z_2$ with $z_1, z_2 \not \in\Pp \cap R_1$ while $z_1, z_2 \in R_1$.  In this case $z_1, z_2 \not \in \Pp$ and we have a contradiction of our assumption that $\Pp$ is a prime ideal of $R_2$.
\end{proof}
\begin{lemma}
\label{lemma:prime}
For any $j =1,\ldots, n$ and any  $a$ we have that $(x_j-a)\bar R_j$ is a prime ideal.
\end{lemma}
\begin{proof}
Observe that for any $a \in k$ we have that $x_j -a$ generates a prime ideal of $k_j[x_j-a]$.  Further, since by assumption $F(z,x_1,x_2,x_{j-1},c,x_j,\ldots,x_n)$ is irreducible  for any $c \in k$, we can also conclude by Proposition 25, Chapter 1, \S 8 of \cite{L}, that $(x_j-a)$ generates a prime ideal in $\bar R_j$, the integral closure of $k_j[x_j-a]$ in $K$.
\end{proof}
From Lemma \ref{le:bigr} and Lemma \ref{lemma:prime} we have the following corollary.
\begin{corollary}
\label{cor:prime}
For any $j=1,\ldots,n$ and any $a \in \phi(I_j)$ we have that $(x_j-a)\bar R_{I_1,\ldots,I_n}$, $ (x_j-a)\bar R, (x_j-a)R$ are all prime ideals.
\end{corollary}
\begin{proof}
By Lemma \ref{le:bigr} it is enough to observe that $k[\bar x] \subseteq R_{I_1,\ldots,I_n} \subseteq k_j(x_j)$, and  $$R \subseteq \bar R \subseteq \bar R_{I_1,\ldots,I_n} \subseteq \bar R_j.$$
\end{proof}
Below are several general results concerning prime ideals.
\begin{lemma}
\label{le:coef}
Let $A$ be a ring and let $\pp$ be a prime ideal of $A$.   If $y \in A \setminus \pp A$ and for some elements of $a_0, \ldots,a_{m-1} \in A$ we have that $y^m+a_{m-1}y^{m-1} +\ldots + a_0=0$, then  some $a_i \not \in \pp A$.
\end{lemma}
\begin{proof}
Assume $y \not \in \pp A$ and $a_0,\ldots,a_{m-1} \in \pp A$ and conclude that $y^m + a_{m-1}y^{m-1} + \ldots + a_1y \in \pp A$.  Since $\pp$ is prime and $y \not \in \pp A$ we conclude that   $y^{m-1} + a_{m-1}y^{m-1} + \ldots + a_1 \in \pp A$ and therefore, since $a_1 \in \pp A$ we have $y^{m-2} + a_{m-1}y^{m-1} + \ldots + a_2 \in \pp A$.  Continuing by induction, we conclude that $y \in \pp A$ contradicting our assumptions on $y$.
\end{proof}
Applying Lemma \ref{le:coef} to our situation we obtain the following.
\begin{lemma}
\label{le:notin}
Let $A, B$ be one of the following pairs of rings: $(\bar R, k[\bar x]), (\bar R_{I_1,\ldots,I_n}, R_{I_1,\ldots,I_n})$ or $(\bar R_{j,a}, R_{j,a})$.     For $y \in A$ we have that $y \in (x_j-a)A$ for some $j =1,\ldots,n$ and some $a \in K$ if and only if all the coefficients, excluding the leading one, of its monic irreducible polynomial over $B$ are in $(x_j-a)B$.
\end{lemma}
\begin{proof}
Let $A^G$ be the integral closure of $A$ in $K^G$.  If $y \in  (x_j-a)A$ and  $y=y_1, \ldots,y_k \in A^G$ are all the conjugates of $y$ over $k(\bar x)$, then $y_i =(x_j-a)u_i$, where $u_1,\ldots, u_k \in A^G$ and are a complete set of conjugates of $u_1$ over $k(\bar x)$.  Thus any symmetric function of $u_i$'s is in $A^G\cap k(\bar x) = B$, and any symmetric function of $y_i$ is in $(x_j-a)B$.  Thus, if $y \in (x_j-a)A$ all the coefficients its monic irreducible polynomial over $B$ are in $(x_j-a)B$.   The converse follows from Lemma \ref{le:coef}.
\end{proof}

Using an approach similar to the one used in Lemma \ref{le:coef}, we can show the following.
\begin{lemma}
\label{le:same}
Let $B \subset A$ be integral domains with $A$ being the integral closure of $B$ in the fraction field of $A$.  If $\pp \subseteq \qq$ are prime ideals of $A$ such that $\qq A \cap B =\pp A \cap B$, then $\pp=\qq$.
\end{lemma}
\begin{proof}
Let $y \in  \qq A \setminus \pp A$ and let $T^m+ a_{m-1}T^{m-1} + \ldots + a_0$ be the monic irreducible polynomial of $y$ over $B$.  Since $y \in \qq A$ we have that $a_0 \in \qq A \cap B \subseteq \pp A$ and therefore $$y^m +a_{m-1}y^{m-1} + \ldots + a_1y \in \pp A.$$  At the same time $y \not \in \pp A$ and therefore we have $y^{m-1} + a_{m-1}y^{m-2} + \ldots + a_1 \in \pp A \subseteq \qq A$. Hence $a_1 \in \qq A \cap B \subset \pp A$.  Continuing by induction we conclude that $y \in \pp A$ contradicting our assumption on $y$.
\end{proof}

We now proof a stronger version of Lemma \ref{le:notin}.
\begin{lemma}
\label{le:norm}
Let $B$ be an integrally closed domain with a fraction field $F$.  Let $H$ be a finite extension of $F$.  Let $A$ be the integral closures of $B$ in $H$.  Let $\pp$ be a prime ideal of $B$ and assume $\pp A$ is also a prime ideal with $\pp A \cap B=\pp B$.  We claim that if $y \in A$, then $y \in \pp A$ if and only if $\bm N_{H/F}y \in \pp B$.
\end{lemma}

\begin{proof}
 Let $H^G$ be the Galois closure of $H$ over $F$ (in some fixed algebraic closure of $F$). Let  $A^G$ be the integral closure of $B$  and $A$ in $H^G$. (This makes sense by Lemma \ref{le:basic}.)  Let $y=y_1, \ldots,y_k$ be all the conjugates of $y$ over $F$ and observe that they all are in $A^G$.  Let $\Jj$ be any prime ideal of $A^G$  such that $\Jj A^G \cap F= \pp B$.  Such a $\Jj$ exists by Lemma \ref{le:basic}. Further, $\Jj A^G \cap A$ contains $\pp B$ and thus contains $\pp A$.   Therefore, by Lemma \ref{le:same}, $\pp A \subseteq \Jj A^G \cap A$  implies that $\Jj A^G \cap A= \pp A$.  For $i=1,\ldots, k$, let $\sigma_i$ be an automorphism sending $y=y_1$ to $y_i$.  Note that $y \in \pp A$ if and only if $y_i \in \pp \sigma_i(A)$ and $\Jj A^G \cap \sigma_i(A) = \pp \sigma_i(A)$ by Lemma \ref{le:same}.  Thus, if $y \not \in \pp A$, for all $i =1,\ldots, k$ we have that $y_i \not \in \pp \sigma_i(A)$ and therefore, $y_i \not \in \Jj A^G$.  Since $\Jj$ is prime, we conclude that in this case  $y_1\ldots y_k \not \in \Jj A^G$ and consequently  $y_1\ldots y_k \not \in \Jj A^G \cap B=\pp B$.
Conversely, if $y \in \pp A$, then $y_i \in \Jj A^G$, for all $i$ and $y_1 \ldots y_k \in \pp B$.
\end{proof}
We also need the following easy fact.
\begin{lemma}
\label{le:power}
Let $A, B$ be one of the following pairs of rings: $(\bar R, k[\bar x]), (\bar R_{I_1,\ldots,I_n}, R_{I_1,\ldots,I_n})$ or $(\bar R_{j,a}, R_{j,a})$, where $j=1,\ldots,n$ and $a \in \phi(I_j)$.  If  $v \in A$ then there exists $k \in \Z_{\geq 0}$ such that
\begin{equation}
\label{eq:exponent}
\frac{v}{(x_j-a)^r} \in A \land \frac{v}{(x_j-a)^{r+1}} \not \in A,
\end{equation}
\end{lemma}
\begin{proof}
Let $y$ and suppose $\displaystyle \frac{y}{(x_j-a)^r} \in A$.  In this case $\displaystyle z=\bm N_{K/k(\bar x)}y \in B \subset k(\bar x)$ and is divisible by $(x_j-a)^{r[K:k(\bar x)]}$ in $B$.  Thus the numerator of $z$, a polynomial, is divisible by $\displaystyle (x_j-a)^r$ and $r$ cannot be arbitrarily large.
\end{proof}
In the future we will call $r$ the order of $y$ at $x_j-a$ and write $r = \ord_{(x_j-a)}y$.  Next we prove another easy property of norms.
\begin{lemma}
\label{le:seminorm}
Let $B$ be an integral domain with a fraction field $F$.  Let $H$ be a finite extension of $F$.  Let $A$ be integral closures of $B$ in $H$. If  $y \in A$ and $N={\bm N}_{H/F} y$, then $\displaystyle \frac{N}{y} \in A$.
\end{lemma}
\begin{proof}
Let $H^G$ be the Galois closure of $H$ over $F$ (in some fixed algebraic closure of $F$). Let  $A^G$ be integral closures of $B$  and $A$ in $H^G$.  Let $y=y_1, \ldots,y_k$ be all the conjugates of $y$ over $F$ and observe that they all are in $A^G$.   Since $N \in B \subset A$ and $y \in A$, we conclude that $\frac{N}{y} \in H$.  At the same time $\frac{N}{y} \in A^G$ and therefore in $A^G \cap H=A$, since $A$ is the set of all elements of $H$ integral over $B$.
\end{proof}
Using Lemma \ref{le:notin} we can also prove that $\bar R_{j,a}$ is a localization. (We remind the reader that $\bar R_{j,a}$ is defined as the integral closure of $R_{j,a}$ in $K$.)
\begin{lemma}
\label{le:local}
$\bar R_{j,a}$ is the localization of $\bar R$ at $x_j-a$.
\end{lemma}
\begin{proof}
First we show that if $y \in \bar R_{j,a}$ it must be of the form $\displaystyle y=\frac{y_1}{y_2}$ where $y_1, y_2 \in \bar R$ and $\ord_{(x_j-a)}y_2 = 0$.    Suppose this assertion is not true.  In this case, by Lemma \ref{le:power} there exists a $y \in \bar R_{j,a}$ of the form $\displaystyle y=\frac{y_1}{y_2}$ where $y_1, y_2 \in \bar R$, $\ell=\ord_{(x_j-a)}y_2 > 0$, and $\ord_{(x_j-a)}y_1 =0$.  Let $y_0=y(x_j -a)^{\ell}$ and note that $y_0 \in \bar R_{j,a} \setminus (x_j-a)\bar R_{j,a}$. Let $T^m+ a_{m-1}T^{m-1} + \cdots + a_0$ be the monic irreducible polynomial of $y$ over $R_{j,a}$ and observe that
 \begin{equation}
 \label{eq:irr}
 T^m+ a_{m-1}(x_j-a)^{\ell}T^{m-1} + \cdots + a_0(x_j-a)^{\ell m}
 \end{equation}
  is the monic a polynomial over $R_{j,a}$ satisfied by $y_0$ .  However, by Lemma \ref{le:notin}, the fact that all the coefficients of \eqref{eq:irr}, except for the leading one, are in $(x_j-a)R_{j,a}$ implies that $y_0 \in (x_j-a)\bar R_{j,a}$, producing a contradiction.

Suppose now $\displaystyle y=\frac{y_1}{y_2}$ where $y_1, y_2 \in \bar R$, $\ord_{(x_j-a)}y_2= 0$.  In this case by Lemma \ref{le:norm}, we have that $N={\bm N}_{K/k(\bar x)}y_2 \not \in (x_j -a)k[\bar  x] \subset (x_j -a)R_{j,a}$.  Observe that $\displaystyle N_0 =\frac{N}{y_2} \in \bar R$  by Lemma \ref{le:seminorm} and therefore $Ny \in \bar R $ by Part 3 of Lemma \ref{le:basic}. Now if
 \begin{equation}
 \label{eq:irr1}
 T^m+ a_{m-1}T^{m-1} + \cdots + a_0
 \end{equation}
is the monic irreducible polynomial of $Ny$ over $\bar R$, then $T^k+ \frac{a_{k-1}}{N}T^{k-1} + \cdots + \frac{a_0}{N^k}$ is the monic irreducible polynomial of $y$ over $k(\bar x)$.  Observe that all the coefficients are in $R_{j,a}$ and thus $y \in \bar R_{j,a}$.
\end{proof}

\begin{proposition}
\label{prop:barspec}
$Z_{I_1,\ldots,I_n}=\{(x_j -a): a \in \phi(I_j)\} \cup \{(0)\}$.
\end{proposition}
\begin{proof}
  First of observe that $\Spec R_{I_1,\ldots,I_n}=\{(x_j -a): a \in \phi(I_j)\}$.  Indeed, let $\Ii \in \Spec R_{I_1,\ldots,I_n}$ and consider $\Ii \cap k[\bar x] \not =\emptyset$.  This is  a prime ideal of $k[\bar x]$ by Lemma \ref{le:bigr}, and it must contain an irreducible element, in our case an irreducible polynomial, which  without loss of generality we can assume to be monic.  However, if this polynomial is not of the form $(x_j -a)$ with $a \in \phi(I_j)$, we conclude that $\Ii$ contains units and thus is equal to the whole ring.  Suppose now that $\Ii \cap k[\bar x]$ contains $x_j - a_1$ and $x_j-a_2$ with $a_1 \not =a_2$.  In this case we conclude as above that $\Ii$ contains units.  Finally suppose $\Ii \cap k[\bar x]$ contains $(x_j-a)$ and $(x_{j'}-a')$ for some indices $j\not = j'$.  In this case $\Ii \cap k[\bar x]$ contains $(x_j-x_{j'}+a'')$, a unit in $R_{I_1,\ldots,I_n}$ and we have a contradiction again.  Thus, $\Ii \cap k[\bar x]$ is generated $(x_j -a)$  for some $j$ and some $a \in \phi(I_j)$.  Next let $z \in \Ii$ and observe that $z =\frac{z_1}{z_2}$, where  $z_1, z_2 \in k[\bar x]$ with $z_2$ not divisible by $(x_j - a)$.  At the same time $z_2z = z_1 \in \Ii \cap k[\bar x]$ and therefore $z_1$ is divisible by $(x_j - a)$.  Therefore, $\Ii$ is generated by $(x_j - a)$ also.

Now let $\Ii \in \Spec \bar R_{I_1,\ldots,I_n}$  and observe that $\Ii \cap k[\bar x] = (x_j-a)k[\bar x]$ for some $j$ and some $a \in \phi(I_j)$, and thus $\Ii\bar  R_{I_1,\ldots,I_n}$ contains $(x_j-a)R_{I_1,\ldots,I_n}$.   At the same time, by Corollary \ref{cor:prime}, $(x_j-a)\bar R_{I_1,\ldots,I_n}$ is also a prime ideal.  Hence, by Lemma \ref{le:same}, we have that $$\Ii\bar R_{I_1,\ldots,I_n}=(x_j-a)\bar R_{I_1,\ldots,I_n}.$$
\end{proof}

\begin{remark}
We also have an alternate description of $R_{I_1,\ldots,I_n}$ and $\bar R_{I_1,\ldots,I_n}$.  We remind the reader that for $j =1,\ldots,n, a \in \phi(I_j)$, it is the case that
 $$R_{a,j}=\{y\in k(\bar x): \frac{y_1}{y_2}, y_1,y_2 \in k[\bar x], y \not \equiv 0 \mod x_j-a\}.$$
Thus, it is not hard to see that
$$R_{I_1,\ldots,I_n} =\bigcap_{a \in \phi(I_j), j=1, \ldots, n} R_{a,j}, $$
$$\bar R_{I_1,\ldots,I_n} =\bigcap_{a \in \phi(I_j), j=1, \ldots, n} \bar R_{a,j}, $$
and by Lemma \ref{le:local} we know that $\bar R_{a,j}$ is the localization of $\bar R$ at $x_j -a$.
\end{remark}

Now that we established the structure of $Z_{I_1,\ldots,I_n}$ we can show the following.
\begin{proposition}
\label{prop:RZ}
$\bar R_{I_1,\ldots,I_n}$ is Turing equivalent to $Z_{I_1,\ldots,I_n}$.
\end{proposition}
\begin{proof}
First of all the fact that $\bar R_{I_1,\ldots,I_n}$ is Turing reducible to $Z_{I_1,\ldots,I_n}$  follows from the fact that  $\bar R_{I_1,\ldots,I_n}$ is the ring corresponding to the whole space.  To establish the reverse reducibility, observe that $\Spec \bar R_{I_1,\ldots,I_n} = \{(x_j -a): \displaystyle \frac{1}{x_j-a} \not \in \bar R_{I_1,\ldots,I_n}\} \cup \{(0)\}$, any closed subset of $Z_{I_1,\ldots,I_n}$ has to be co-finite, since any element of the ring contained in the intersection of a collection of ideals of the form  $(x_j -a)$, must be divisible by all $x_j-a$, and thus the collection must be finite.  Finally if $J_1, \ldots, J_n$  are finite, possibly empty, subsets of positive integers and  $(I_1', \ldots,I'_n) = (I_1\setminus J_1,\ldots,I_n\setminus J_n)$, then the ring corresponding to $Z_{I_1,\ldots,I_n} \setminus \{(x_j-a), a \in \phi(J_j)\}$ is $\bar R_{I_1', \ldots,I'_n}$, computable uniformly in $(J_1, \ldots,J_n)$ from  $\bar R_{I_1,\ldots,I_n}$, and thus from $\bar R_{I_1,\ldots,I_n}$ only, since $(J_1, \ldots,J_n)$ is finite.
\end{proof}
As above, the argument of Proposition \ref{prop:RZ} can be easily adapted to show the enumerable equivalence as well.
\begin{lemma}
$\bar R_{I_1,\ldots,I_n}$ is enumeration equivalent to $Z_{I_1,\ldots,I_n}$.
\end{lemma}
Now from Proposition \ref{prop:Turing1} we get the following corollary.
\begin{corollary}
\label{cor:Turing}
$Z_{I_1,\ldots,I_n}$ is Turing and enumeration reducible to  $\max(I_1,\ldots,I_n)$.
\end{corollary}

To state the final results of this section, we need additional notation:
\begin{notation}\rule{5pt}{0in}
\begin{itemize}
\item For $I \subset \Z_{>0}$, let $R_I =R_{I_1,\ldots,I_n}$, where $I_1 = I$  and $I_j = \Z$ for all $j=2,\ldots,n$.
\item Define $\bar R_I, Z_I$ in a similar fashion.
\item Let $Z=\Spec \bar R_{\Z}$.
\end{itemize}
\end{notation}
From  Corollary \ref{cor:Turing} we now derive the following result.
\begin{theorem}
\label{thm:scheme1} 
The class of isomorphism types of unions of
subspaces of $(Z,O_Z)$ satisfies conditions of Proposition \ref{prop:schemes} and thus admits arbitrary degrees of isomorphism types,
and admits isomorphism types without degree.
\end{theorem}

We now describe some examples of $V$ satisfying our conditions (besides $V\cong \A_k^n$).  First we need a simple lemma.
\begin{lemma}
\label{le:ex}
If $G(X)$ is a rational function field over a field of constants $G$,  $\alpha$ a root of a polynomial
\begin{equation}
\label{eq:example}
Q(T) -P(X)=0
\end{equation}
in some algebraic closure of $G(X)$, where $P(X) \in G[X], Q(T) \in G(T), \deg Q(T)=m$ and $(\deg P(X),m)=1$, then $[G(\alpha, X):G(X)]=m$.
\end{lemma}
\begin{proof}
Observe that $\alpha$ is integral over $G[X]$ and therefore it can have a pole only at a $G(\alpha,X)$-valuation $\pp_{\infty}$  lying above the infinite valuation of $G(X)$ (otherwise known as the ``degree'').   Thus $\ord_{\pp_{\infty}}Q(\alpha)= m\ord_{\pp_{\infty}}\alpha = e\deg P(X)$, where $e$ is the ramification degree of $\pp_{\infty}$ over the infinite valuation of the rational function field.  Since $(m,  \deg P(X))=1$, we must conclude that $m \mid e$ and therefore the degree of the extension is at least $m$.  However, $\alpha$ is of degree at most $m$ over $G(X)$.  Thus the degree of the extension is exactly $m$.
\end{proof}
\begin{example}
Let  $F(z,x_1,\ldots,x_n) \in k[x_1,\ldots, x_n, z]$  satisfy the following conditions.
 \begin{itemize}
 \item $F(z,x_1,\ldots,x_n)=Q(z)-A(x_1,\ldots, x_n)$, $Q(z) \in k[z], A(x_1,\ldots, x_n) \in k[x_1,\ldots,x_n]$
 \item $Q(z)$ is of degree $m$.
 \item Considered as a polynomial in $x_j, j=1,\ldots, n$,  the polynomial $A(x_1,\ldots,x_n)$ is monic and of degree $d_j$.
 \item $(d_j,m)=1$.
 \end{itemize}
   Let $n \in \Z_{>1}$,  $r,j=1,\ldots,n, r \not=j$, $G=k(x_1,\ldots,x_{r-1},x_r,\ldots,x_n)$. Note that $k(\bar x)=G(x_r)$ and $F(z,x_1,\ldots,x_n)=Q(z)-P(x_r)$, where $Q(z)\in G[z], P(x_r)\in G[x_r]$, the second polynomial is monic, and $(\deg P, \deg Q)=1$.  Substituting an element  of $k$ for $x_j, j\not=r$ does not change the degree of $Q$ in $z$ or the degree of $P$ in $x_r$ since $P$ is monic.  Hence  $F(z,x_1,\ldots,x_{j-1},a,x_j,\ldots,x_n)$ is of the form \eqref{eq:example} and  we can apply Lemma \ref{le:ex} to conclude that $$[k(z_j(a),x_1,\ldots,x_{j-1},x_j,\ldots,x_n):k(x_1,\ldots,x_{j-1},x_j,\ldots,x_n)]=m,$$ where  $z_j(a)$ is a root of $F(z,x_1,\ldots,x_{j-1},a,x_{j+1},\ldots,x_n)$ in the algebraic closure of $$k(x_1,\ldots,x_{j-1},x_j,\ldots,x_n).$$  Thus, the irreducibility condition in Notation and Assumptions \ref{not:irr} is satisfied.
\end{example}

 We finish with somewhat different example of the  ``subspace'' construction.
 \begin{example}
 Let $Z = \Z_{>0}$ under the co-finite topology.  Let  $R(Z)$ be the ring of all rational functions over $\Q$ defined on  $\Z>0$.  Let $Z_X=\{ 2n+1| n \in X\} \cup \{2n| n \in \Z_{n>0}\}$.    Using the same argument as in Proposition \ref{prop:schemes}, we can show that  $(Z,O_Z)$ together with the collection $\{(Z_X, O_{Z_X})\}, X \subset \Z_{>0}$ satisfies the assumptions of Theorem \ref{thm:subspace}.\\

 We can generate more examples of this sort.  All we need is a countably infinite domain we can convert into a topological space using co-finite topology and a computable  field of functions, where one can determine the domain of a function under a uniform effective procedure.
 \end{example}

\appendix

 \section{A Family of Non-Isogenous Curves}%
\label{sec:ecurves}%
\setcounter{equation}{0}
  In this section we construct a computable collection of pairwise
  non-isogenous curves, which we use in Example \ref{needapp1}.

\begin{lemma}%
\label{le:primes}%
Let $K$ be a number field. Then there exists a computable sequence of its primes such that the equation $Z^3-4=0$.
has solutions in the residue field of each prime in the sequence. %
\end{lemma}%
\begin{proof}%
If a prime $\mathsf{p}$ of $K$ is not a factor of 2 or 3, then $Z^{3}-4=0$ has solutions mod
$\mathsf{p}$ if and only if $\mathsf{p}$ has a relative degree one factor in the extension
$K(\sqrt[3]{4})/K$. In particular, $Z^{3}-4=0$ has solutions mod $\mathsf{p}$ if the prime in
question splits completely in the extension $K(\sqrt[3]{4})/K$. (See Proposition 25, page 27 of
\cite{L}.) By the Chebotarev Density Theorem (see Theorem 10.4, page 182 of \cite{Januz}), there are
infinitely many such primes. Using a computable listing of all primes of $K$ as in \cite{CHS} and a
computable listing of algebraic integers, we can generate a computable listing of primes of $K$
which have $\sqrt[3]{4}$ in their residue field, using the fact that evaluating the order at a
given prime is a computable operation. (See \cite{CHS}.)
\end{proof}%
\begin{lemma}%
\label{le:sequence}%
Let $\{{}\mathsf{p}_{l}\}$ be a computable listing of $K$ primes excluding primes dividing 2 and 3 and such that
the residue field of every prime in the sequence contains $\sqrt[3]{4}$. Then there exists a computable sequence
$\{A_{i}\}\in O_{K}$ satisfying the following computable conditions. %
\begin{enumerate}%
\item $4A_1^3+27\not=0,4A_1^3+27\cong 0\mod\mathsf{p} _{1} $. %
\item Let $\mathcal{G}(k)$ be the set of primes of $K$ dividing $4A_i^3+27$ for $i=1,\ldots ,k-1$. Then for any
prime $\mathsf{p}\in \mathcal{G}(k)$,we have that $\ord_{\mathsf{p}}(4A_k^3+27)=0$. Moreover, if  $\mathsf{q}$ is
the first prime in the listing above which is not in $\mathcal{G}(k)$. Then $\ord_{\mathsf{q}}(4A_k^3+27)>0$.%
\end{enumerate}%
\end{lemma}%

\begin{proof}%
Let $A_{1}$ be an algebraic integer of $O_{K}$ such that $A_{1}\equiv \frac{-3}{\sqrt[3]{4}} \mod \mathsf{p}_1$
but $A_{1}\not=\frac{-3}{\sqrt[3]{4}}$. Such a $A_{1}$ exists by the Strong Approximation Theorem and clearly
satisfies the requirements. Next assume $A_{1},\ldots ,A_{k-1}$ have been defined. Let
$\mathcal{H}(k)=\prod_{\mathsf{p}\in \mathcal{G}(k)} \mathsf{p}$. Then, by the Strong Approximation Theorem
again, there exists an algebraic integer $A_{k}$ such that $A_{k}\equiv 0\mod\mathcal{H}(k)$,
$A_{k}\not=\frac{-3}{\sqrt[3]{4}}$, and $A_{k}\equiv \frac{-3}{\sqrt[3]{4}} \mod~\mathsf{p}_{k}$. Since the
requisite numbers exist, we can find them by a systematic search of $O_{K}$ provided we can check that a number
satisfies the requirement in an effective manner. The only operation which can be a source of difficulty is
evaluation of order at a prime. However, we have addressed that issue in \cite{CHS}.%
\end{proof}%
\begin{lemma}%
\label{le:isogeny}%
Let $\{A_i\}$ be a sequence of $O_K$-integers defined in Lemma \ref{le:sequence}. Let $\{E_i\}$ be a sequence of
curves defined over $ K$, such that the affine part of $E_i$ is given by equation $y^2=x^3+A_ix +1$ . Then the
sequence $E_i$ satisfies the following conditions.%
\begin{enumerate}%
\item \label{it:el1}%
For all $i$ we have that $E_i$ is an elliptic curve.%
\item \label{it:el2} The equation $y^{2}=x^{3}+A_ix+1$ is minimal with respect to all primes not dividing $2$
or $3$.%
\item \label{it:el3} For every $i\in \omega$ there exists a $K$-prime $ \mathsf{q}_{i}$ such that $E_{i}$ does not
have a potentially good reduction at $\mathsf{q}_{i}$, but for any $k\in \omega, k\not=i$, we have that $E_{i}$
has a good reduction at $\mathsf{q}_{k}$.%
\item \label{it:el4} For any pair of natural numbers $i \not = k$, we have $ E_i$ and $E_k$ are
not isogenous over $\tilde{K}$ -- algebraic closure of $K$.%
\end{enumerate}%
\end{lemma}%
\begin{proof}%
(\ref{it:el1}) All the curves in our sequence are defined by Weierstrass equations with non-zero discriminants. (For a
Weierstrass equation of the form $y^2=x^3+Ax+B$, the discriminant $\triangle=-16(4A^3+27B^2)$. See Remark 1.3, pages 49--50 of
\cite{Sil1}.) Thus, for every $i$ we have that $E_i$ is an elliptic curve by Proposition 1.4(a), page 50 and Proposition 3.1,
page 63 of \cite{Sil1}. \\%
(\ref{it:el2}) This part follows from Remark 1.1, page 172 of \cite{Sil1} if we remind the reader that given the form of our
Weierstrass equations, $ c_6=-216$. \\%
(\ref{it:el3}) First of all, we observe that ${}\mathsf{q}_{i}$ divides $ \triangle (E_{i})$ by construction. Since the equation
we have constructed for $E_{i}$ is minimal at ${}\mathsf{q}_{i}$, by Proposition 5.1, page 180 of \cite{Sil1}, we know that
$E_{i}$ does not have a good reduction at $ {\mathsf q}_{i}$. Further, $j(E_{i})=\frac{-1728(4^{3})}{\triangle (E_{i})}$ is not
integral at ${}\mathsf{q}_{i}$. Then by Proposition 5.5, page 181 of \cite {Sil1}, $E_i$  does not
have a potential good reduction at ${\mathsf q}_{i}$. Finally, if $k\not=i$, then $\triangle (E_{k})\not\equiv
0~\mod~\mathsf{q}_{i}$ and by Proposition 5.1, page 180 of \cite{Sil1} again, $E_{k}$ does have a good reduction modulo
${\mathsf q}_{i}$\\%
(\ref{it:el4}) If $E_i$ is isogenous to $E_k$ over some number field $ K^{\prime}$ then they have good reduction at the same
primes of $K^{\prime}$. (Corollary 7.2, page 185 of \cite{Sil1}.) On the other hand, by construction and definition of
potential good reduction, for any $K^{\prime}$ -- finite extension of $K$, for $i\not = k$, we have that $E_i$ and $E_k$ will
have distinct sets of good primes (and bad primes). Therefore, $E_i$ and $E_k$ will not be isogenous over any finite extension
of $K$. Thus they will not be isogenous over
$\tilde{K}$.%
\end{proof}%

\end{document}